\documentclass[pdflatex,sn-mathphys-num]{sn-jnl}

\usepackage{amsmath}
\usepackage{amssymb}
\usepackage{enumerate}
\usepackage{algorithm,algorithmic}
\usepackage{caption}
\usepackage{array}
\usepackage{pgfplots, pgfplotstable, booktabs, colortbl, array}
\usepackage{multirow}

\usepackage{float}
\floatstyle{plaintop}
\restylefloat{table}

\usepackage{todonotes}

\newcommand{\ignore}[1]{}

\usepackage{amsfonts}
\usepackage{amssymb}

\newcommand{\norm}[1]{\|#1\|}
\newcommand{\abs}[1]{|#1|}

\newcommand{\uu}{{\bf u}}


\newtheorem{theorem}{Theorem}[section]
\newtheorem{lemma}[theorem]{Lemma}

\theoremstyle{definition}
\newtheorem{definition}[theorem]{Definition}

\theoremstyle{remark}

\theoremstyle{thmstyleone}%

\theoremstyle{thmstyletwo}%

\theoremstyle{thmstylethree}%

\begin{document}

\title[RCLUPP]{A new randomized CholeskyQR based on LU decomposition with partial pivoting}


\author[1]{\fnm{Yuwei} 
\sur{Fan}}\email{fanyuwei2@huawei.com}

\author*[1]{\fnm{Haoran} \sur{Guan}}\email{guan.haoran@huawei.com}

\author[2]{\fnm{Zhonghua} \sur{Qiao}}\email{zhonghua.qiao@polyu.edu.hk}

\affil[1]{\orgdiv{Theory Lab}, \orgname{Huawei Leibniz Research Center}, \orgaddress{\street{Sha Tin}, \city{New Territories}, \postcode{999077}, \state{Hong Kong SAR}, \country{China}}}

\affil[2]{\orgdiv{Department of Applied Mathematics}, \orgname{The Hong Kong Polytechnic University}, \orgaddress{\street{Hung Hom}, \city{Kowloon}, \postcode{999077}, \state{Hong Kong SAR}, \country{China}}}


\abstract{CholeskyQR has received considerable attention in recent years for its efficiency and simplicity in computing QR decomposition of the tall-skinny $X \in \mathbb{R}^{m\times n}$ with $m \ge n$ and $\mbox{rank}(X)=n$. Leveraging matrix sketching from randomized linear algebra, randomized CholeskyQR (RCholeskyQR) has been proposed to accelerate the computation by reducing the dimension of the problems. In this work, we propose RCLUPP, a new randomized CholeskyQR-type  algorithm based on LU decomposition with partial pivoting (LUP decomposition). By taking LUP decomposition and the thin HouseholderQR on the sketched matrix, RCLUPP significantly improves the applicability and efficiency compared with LU-CholeskyQR2 (LC2). We present a rigorous rounding error analysis of RCLUPP, with a sharper bound of residual compared with those in the existing works. Comparative studies demonstrate that RCLUPP outperforms CholeskyQR2, Shifted CholeskyQR3 (SCholeskyQR3), and LC2 in terms of applicability while maintaining competitive accuracy and efficiency. A variant, RCLUPPr, performs LUP decomposition directly on $X \in \mathbb{R}^{m\times n}$, offering exceptional robustness and numerical stability for the ill-conditioned scenarios, which exceeds that of RCLUPP and RCholeskyQR. Numerical experiments on the synthetic and real-world matrices validate the theoretical results.}




\maketitle

\section{Introduction}
QR decomposition is a fundamental tool in numerical linear algebra with broad applications. This paper focuses on QR decomposition of the tall-skinny $X \in \mathbb{R}^{m\times n}$ with $m \ge n$ and $\mbox{rank}(X)=n$ in the form of
\begin{equation}
X = QR, \nonumber
\end{equation}
where $R \in \mathbb{R}^{n\times n}$ is upper-triangular and $Q \in \mathbb{R}^{m\times n}$ is orthogonal and numerically well-conditioned.

A variety of algorithms for QR decomposition exist, including HouseholderQR, classical Gram--Schmidt (CGS), modified Gram--Schmidt (MGS) and TSQR \cite{ballard2011, 2011, Communication, MatrixC, Higham, Numerical}. CholeskyQR offers advantages over these methods through its reliance on BLAS3 operations and its simple reductions in a parallel environment. In particular, CholeskyQR-type algorithms are communication-avoiding (CA), requiring only one or two global reductions regardless of the matrix height, which makes them highly attractive for large-scale parallel scientific computing. We begin with the basic CholeskyQR, as outlined in Algorithm~\ref{alg:cholqr}. For $X \in \mathbb{R}^{m\times n}$ with $m \ge n$, the upper-triangular factor $R$ is obtained from Cholesky decomposition of the gram matrix $X^{\top}X$. The orthogonal factor $Q \in \mathbb{R}^{m\times n}$ is then computed by solving linear systems, which can be written $Q=XR^{-1}$.

\begin{algorithm}
\caption{$[Q,R]=\mbox{CholeskyQR}(X)$}
\label{alg:cholqr}
\begin{algorithmic}[1]
\REQUIRE $X \in \mathbb{R}^{m\times n}$ with $\mbox{rank}(X)=n$.
\ENSURE \mbox{Orthogonal factor} $Q \in \mathbb{R}^{m\times n}$, \mbox{Upper triangular factor} $R \in \mathbb{R}^{n \times n}.$ 
\STATE $G=X^{\top}X,$
\STATE $R=\mbox{Cholesky}(G),$
\STATE $Q=XR^{-1}.$
\end{algorithmic}
\end{algorithm}%

The applicability and accuracy of a single CholeskyQR is limited in the real implementation. For $X$ in a moderate size, the orthogonality deteriorates rapidly when $\kappa_{2}(X)$ approaches $\uu^{-\frac{1}{2}}$. CholeskyQR2 \cite{2014, error} improves the accuracy by repeating the procedure twice, yet it still suffers from the potential numerical breakdown when $\kappa_{2}(X) \ge \uu^{-\frac{1}{2}}$. To address this problem, Shifted CholeskyQR3 (SCholeskyQR3) \cite{Shifted} has been developed with a shifted item $s$ in Cholesky decomposition, but it remains ineffective for many ill-conditioned cases due to the choice of $s$ and the imposed bound on $\kappa_{2}(X)$. More recently, with the rapid development of randomized linear algebra \cite{halko2011, 2020}, randomized CholeskyQR (RCholeskyQR) has been proposed and analyzed in \cite{Randomized, Novel, Householder}, which integrates matrix sketching with different strategies to compute the $R$-factor. Although these methods achieve a better balance among the applicability, accuracy and efficiency, they may still fail when $\kappa_{2}(X)$ is sufficiently large, consistent with their sufficient conditions of $\kappa_{2}(X)$ theoretically.

Another line of the development of CholeskyQR is LU-CholeskyQR2 (LC2) \cite{LUChol}, which introduces LU decomposition with partial pivoting (referred to simply as LUP decomposition) as a preconditioning step. The $R$-factor is then obtained from the upper-triangular factor of LUP decomposition combined with a subsequent step of Cholesky decomposition. Compared with other CholeskyQR-type algorithms, LC2 handles some ill-conditioned matrices more effectively and achieves superior accuracy. However, it still fails on the matrices with some specific structures and incurs additional computational cost due to LUP decomposition of $X \in \mathbb{R}^{m\times n}$ directly.

Several strategies can be employed to mitigate the limitations of LC2. The randomized technique, such as matrix sketching \cite{rgs}, accelerates LUP decomposition by reducing the size of the problem. In addition, some alternative methods for computing the upper-triangular factor, rather than Cholesky decomposition of the gram matrix, can avoid numerical breakdown for the ill-conditioned matrices or those with special structures. Consequently, the resulting randomized algorithm achieves superior applicability compared with LC2, while retaining comparable accuracy and efficiency with matrix sketching.

In this work, we propose RCLUPP, a new randomized CholeskyQR-type algorithm based on LUP decomposition for the tall skinny $X$ with $m \ge n$ and $\mbox{rank}(X)=n$, as shown in Algorithm~\ref{alg:RCLUPP}. In contrast to LC2, we apply matrix sketching prior to LUP decomposition to reduce the dimension of the problem and accelerate the algorithm. Moreover, Cholesky decomposition is replaced by a more stable way for computing the R-factor alternatively, such as the thin HouseholderQR or Givens rotations. We focus on the case with HouseholderQR in this work. After the preconditioning step, CholeskyQR2 is applied to guarantee the accuracy and numerical stability of the algorithm. By compressing the input $X \in \mathbb{R}^{m\times n}$ into a smaller $X_{s} \in \mathbb{R}^{s\times n}$ with $s=\mathcal{O}(n)$, both LUP decomposition and the thin Householder QR of $X_{s}$ become significantly more efficient. Although LUP decomposition is introduced, it is performed exclusively on $X_{s} \in \mathbb{R}^{s\times n}$ with $s=\mathcal{O}(n)$, thereby adding negligible communication cost to the original $X \in \mathbb{R}^{m\times n}$ and preserving the CA nature of CholeskyQR. We present a detailed rounding error analysis demonstrating the numerical stability of RCLUPP, together with a tighter bound of residual compared with those in the existing works and a discussion of the sufficient condition for the applicability. A variant of RCLUPP, named RCLUPPr, performs LUP decomposition directly on the original $X$ to eliminate the effect of matrix sketching on the applicability. We have a recent subsequent work, \cite{RCLUPPr}, providing a detailed analysis of RCLUPPr and its implementation. Numerical experiments show that RCLUPP achieves better applicability than CholeskyQR2, SCholeskyQR3 and LC2, together with superior accuracy compared with that of RCholeskyQR in many practical cases. When $\frac{m}{n}$ is large and $s<<m$, RCLUPP is faster than LC2 while maintaining the CPU time (s) comparable to those of SCholeskyQR3 and RCholeskyQR with appropriate optimizations. For the very ill-conditioned cases, RCLUPPr, though slightly slower, delivers outstanding applicability and numerical stability. In summary, we propose a new family of the randomized CholeskyQR-type algorithms that incorporate LUP decomposition. These algorithms achieve a balance among the applicability, accuracy and efficiency for QR decomposition of the tall-skinny $X \in \mathbb{R}^{m\times n}$ with $m \ge n$ and $\mbox{rank}(X)=n$.

\begin{algorithm}
\caption{$[Q,R]=\mbox{RCLUPP}(X)$}
\label{alg:RCLUPP}
\begin{algorithmic}[1]
\REQUIRE $X \in \mathbb{R}^{m\times n}$ with $\mbox{rank}(X)=n$, \mbox{Sketch matrix} $S \in \mathbb{R}^{s\times m}$.
\ENSURE \mbox{Orthogonal factor} $Q \in \mathbb{R}^{m\times n}$, \mbox{Upper triangular factor} $R \in \mathbb{R}^{n \times n}.$ 
\STATE $X_{s}=SX,$
\STATE $[M,N,P]=\mbox{LU}(X_{s}),$
\STATE $[Q_{s},R_{0}]=\mbox{HouseholderQR}(L),$
\STATE $R_{1}=R_{0}N,$
\STATE $Q_{1}=XR_{1}^{-1},$
\STATE $[Q,Z]=\mbox{CholeskyQR2}(Q_{1}),$
\STATE $R=ZR_{1}.$
\end{algorithmic}
\end{algorithm}%

This work is the first to develop an efficient RCLUPP algorithm by integrating matrix sketching, LUPP factorization and the thin Householder QR factorization. It also establishes an improved bound of residual and discusses the sufficient condition for the applicability of the algorithm. As a subsequent and complementary work, \cite{RCLUPPr} investigates RCLUPPr from the perspective of the mixed precision with a detailed theoretical analysis, demonstrating its robustness for the ill-conditioned problems and proposing several strategies for the acceleration. Thus, the two works are complementary in the scope while remaining independent in the content.

This paper is organized as follows. In Section~\ref{sec:development}, we review the development of CholeskyQR. Section~\ref{sec:RCLUPP} introduces the proposed RCLUPP and presents a detailed rounding error analysis of the algorithm. Comparisons of the theoretical results and some extensions are shown in Section~\ref{sec:comparison}. Numerical experiments demonstrating the theoretical results are reported in Section~\ref{sec:numerical}. Finally, concluding remarks and discussions are given in Section~\ref{sec:conclusions}.

We use the following notation throughout. For $X \in \mathbb{R}^{m\times n}$, $\norm{X}_{F}$ and $\norm{X}_{2}$ denote the Frobenius norm and spectral norm, respectively. Assume $m\ge n$ and $\mbox{rank}(X)=n$. The singular values of $X$ are ordered as
$\sigma_1(X) \ge \cdots \ge \sigma_n(X)$, where $\sigma_1(X)=\norm{X}_{2}$. The spectral condition number is defined by $\kappa_2(X)=\frac{\sigma_1(X)}{\sigma_n(X)}$. All computations are carried out in the double precision floating-point arithmetic with the unit roundoff $\uu=2^{-53}$. For a matrix $X$, $\abs{X}$ denotes its elementwise absolute value and $fl(\cdot)$ denotes the result of a floating-point computation. The identity matrix of the order $n$ is denoted by $I_n \in \mathbb{R}^{n\times n}$. Finally, $\hat{A}$ denotes the computed counterpart of $A$ in an actual implementation.

\section{Development of CholeskyQR}
\label{sec:development}
In this part, we present the recent development of CholeskyQR with some typical algorithms.

\subsection{CholeskyQR2}
CholeskyQR2 \cite{2014,error} improves the numerical stability of CholeskyQR, particularly its orthogonality, by repeating the procedure twice, as outlined in Algorithm~\ref{alg:cholqr2}. Here, $Q_{1}, Q \in \mathbb{R}^{m\times n}$ denote the $Q$-factors obtained after the first and second steps, respectively, while $R_{1}, Z \in \mathbb{R}^{n\times n}$ are the corresponding $R$-factors.

\begin{algorithm}
\caption{$[Q,R]=\mbox{CholeskyQR2}(X)$}
\label{alg:cholqr2}
\begin{algorithmic}[1]
\REQUIRE $X \in \mathbb{R}^{m\times n}$ with $\mbox{rank}(X)=n$.
\ENSURE \mbox{Orthogonal factor} $Q \in \mathbb{R}^{m\times n}$, \mbox{Upper triangular factor} $R \in \mathbb{R}^{n \times n}.$ 
\STATE $[Q_{1},R_{1}]=\mbox{CholeskyQR}(X),$
\STATE $[Q,Z]=\mbox{CholeskyQR}(Q_{1}),$
\STATE $R=ZR_{1}.$
\end{algorithmic}
\end{algorithm}%

For $X \in \mathbb{R}^{m\times n}$, the orthogonality of a single CholeskyQR can be estimated as
\begin{equation}
\norm{\hat{Q}^{\top}\hat{Q}-I_{n}}_{F} \le \frac{5}{64}\delta^{2}, \nonumber
\end{equation}
with the computed orthogonal factor $\hat{Q} \in \mathbb{R}^{m\times n}$. Here, we have $\delta=8\kappa_{2}(X)\sqrt{mn\uu+n(n+1)\uu}$. For CholeskyQR2, theoretical results of detailed rounding error analysis are shown in Lemma~\ref{lemma 3.1} from \cite{error}.

\begin{lemma}[Rounding error analysis of CholeskyQR2]
\label{lemma 3.1}
For $X \in \mathbb{R}^{m\times n}$ and $[\hat{Q},\hat{R}]=\mbox{CholeskyQR2}(X)$, with
$8\kappa_{2}(X)\sqrt{mn\uu+n(n+1)\uu} \le 1$, $mn\uu \le \frac{1}{64}$ and $n(n+1)\uu \le \frac{1}{64}$, we have
\begin{align}
\norm{\hat{Q}^{\top}\hat{Q}-I_{n}}_{F} &\le 6(mn\uu+n(n+1)\uu), \nonumber \\
\norm{\hat{Q}\hat{R}-X}_{F} &\le 5n^{2}\sqrt{n}\uu\norm{X}_{2}. \nonumber
\end{align}
\end{lemma}

\subsection{SCholeskyQR3}
Although CholeskyQR2 improves the orthogonality by repeating the procedure twice, it still suffers from the numerical breakdown for the ill-conditioned $X$ satisfying $\kappa_{2}(X) \ge \uu^{-\frac{1}{2}}$. To address this limitation, SCholeskyQR introduces a shift item $s$ in Cholesky decomposition to prevent the numerical breakdown \cite{Shifted}. For the accuracy, SCholeskyQR3 combines SCholeskyQR with a subsequent CholeskyQR2, as shown in Algorithm~\ref{alg:Shifted3}.

\begin{algorithm}
\caption{$[Q,R]=\mbox{SCholeskyQR3}(X)$}
\label{alg:Shifted3}
\begin{algorithmic}[1]
\REQUIRE $X \in \mathbb{R}^{m\times n}$ with $\mbox{rank}(X)=n$.
\ENSURE \mbox{Orthogonal factor} $Q \in \mathbb{R}^{m\times n}$, \mbox{Upper triangular factor} $R \in \mathbb{R}^{n \times n}.$ 
\STATE $G=X^{\top}X,$
\STATE take $s>0$,
\STATE $R_{1}=\mbox{Cholesky}(G+sI_{n}),$
\STATE $Q_{1}=XR_{1}^{-1},$
\STATE $[Q,Z]=\mbox{CholeskyQR2}(Q_{1}),$
\STATE $R=ZR_{1}.$
\end{algorithmic}
\end{algorithm}%

The performance of SCholeskyQR3 depends critically on the choice of the shifted item $s$, which is derived from the rounding error analysis of the algorithm. Following the improved bounds in \cite{Columns} than those in \cite{Shifted}, we adopt $s=11(mn\uu+n(n+1)\uu)\norm{X}_{g}^{2}$,  where $\norm{X}_{g}$ is defined as

\begin{definition}
\label{def:g}
For $X=[X_{1},X_{2}, \cdots X_{n-1},X_{n}]\in R^{m\times n}$,
\begin{equation} 
\begin{split}
\norm{X}_{g}:=\max_{1 \le j \le n}\norm{X_{j}}_{2}, \nonumber
\end{split}
\end{equation}
with
\begin{equation}
\begin{split}
\norm{X_{j}}_{2}=\sqrt{x_{1,j}^{2}+x_{2,j}^{2}+……+x_{m-1,j}^{2}+x_{m,j}^{2}}. \nonumber
\end{split}
\end{equation}
\end{definition}

Our improved shift item $s$, which provides better applicability than the original choice $s=11(mn\uu+n(n+1)\uu)\norm{X}_{2}^2$, leads to the following theoretical results for SCholeskyQR3, as stated in Lemma~\ref{lemma 3.2} of \cite{Columns}.

\begin{lemma}[Rounding error analysis of SCholeskyQR3 with the improved $s$]
\label{lemma 3.2}
For $X \in \mathbb{R}^{m\times n}$ and $[\hat{Q},\hat{R}]=\mbox{SCholeskyQR3}(X)$, with $s=11(mn\uu+n(n+1)\uu)\norm{X}_{g}^{2}$, $mn\uu \le \frac{1}{64}$, $n(n+1)\uu \le \frac{1}{64}$ and $\kappa_{2}(X)=\frac{1}{86p(mn\uu+n(n+1)\uu)}$, we have
\begin{align}
\norm{\hat{Q}^{\top}\hat{Q}-I_{n}}_{F} &\le 6(mn\uu+n(n+1)\uu), \nonumber \\
\norm{\hat{Q}\hat{R}-X}_{F} &\le 6.57n^{2}\uu\norm{X}_{g}+4.81n^{2}\uu\norm{X}_{2}. \nonumber
\end{align}
Here, we let $p=\frac{\norm{X}_{g}}{\norm{X}_{2}}$ with $\frac{1}{\sqrt{n}} \le p \le 1$.
\end{lemma}

\subsection{LC2}
Although SCholeskyQR3 improves the applicability over its predecessors, it still fails for many ill-conditioned cases. An alternative approach to enhance the applicability is LU-CholeskyQR2 (LC2) \cite{LUChol}, which incorporates LUP decomposition as a preconditioning step, as outlined in Algorithm~\ref{alg:LC2}. Here, $P \in \mathbb{R}^{m\times m}$ is the permutation matrix, while $M \in \mathbb{R}^{m\times n}$ and $N \in \mathbb{R}^{n\times n}$ are the lower-triangular and upper-triangular factors, respectively.

\begin{algorithm}
\caption{$[Q,R]=\mbox{LC2}(X)$}
\label{alg:LC2}
\begin{algorithmic}[1]
\REQUIRE $X \in \mathbb{R}^{m\times n}$ with $\mbox{rank}(X)=n$. 
\ENSURE \mbox{Orthogonal factor} $Q \in \mathbb{R}^{m\times n}$, \mbox{Upper triangular factor} $R \in \mathbb{R}^{n \times n}.$ 
\STATE $[M,N,P]=\mbox{LU}(X),$
\STATE $G=M^{\top}M,$
\STATE $S=\mbox{Cholesky}(G),$
\STATE $R_{1}=SN,$
\STATE $Q_{1}=XR_{1}^{-1},$
\STATE $[Q,Z]=\mbox{CholeskyQR}(Q_{1}),$
\STATE $R=ZR_{1}.$
\end{algorithmic}
\end{algorithm}%

The theoretical results are shown in Lemma~\ref{lemma 3.3} from \cite{LUChol}.

\begin{lemma}[Rounding error analysis of LC2]
\label{lemma 3.3}
For $X \in \mathbb{R}^{m\times n}$ and $[\hat{Q},\hat{R}]=\mbox{LC2}(X)$, with $mn\uu \le \frac{1}{64}$, $n(n+1)\uu \le \frac{1}{64}$, $8\kappa_{2}(\hat{M})\sqrt{mn\uu+n(n+1)\uu} \le 1$ and $64n^{2}\uu \cdot \kappa_{2}(\hat{M})\kappa_{2}(\hat{N}) \le 1$, we have
\begin{align}
\norm{\hat{Q}^{\top}\hat{Q}-I_{n}}_{F} &\le 6.5(mn\uu+n(n+1)\uu), \nonumber \\
\norm{\hat{Q}\hat{R}-X}_{F} &\le 4.09n^{2}\uu\norm{X}_{2}. \nonumber
\end{align}
Here, $\hat{M}$ and $\hat{N}$ are the computed $M$ and $N$ from LUP decomposition.
\end{lemma}

\subsection{RCholeskyQR}
In recent years, randomized techniques, particularly matrix sketching, have been successfully applied to QR decomposition. Fan et al. \cite{Novel} introduces randomized CholeskyQR (RCholeskyQR) for the first time, while Balabanov \cite{Randomized} provides a detailed analysis of this topic. Building upon these works, Higgins et al. \cite{Householder} utilizes a multi-sketching strategy consisting of two consecutive steps of matrix sketching on CholeskyQR. The general form of RCholeskyQR using a single sketch followed by a thin HouseholderQR is presented in Algorithm~\ref{alg:RC}. Alternatively, Givens rotations can be employed to compute the $R$-factor, $R_{1}$, in the preconditioning step. In this setting, $Q_{s} \in \mathbb{R}^{s\times n}$ denotes the $Q$-factor obtained from the thin HouseholderQR.

\begin{algorithm}
\caption{$[Q,R]=\mbox{RCholeskyQR}(X)$}
\label{alg:RC}
\begin{algorithmic}[1]
\REQUIRE $X \in \mathbb{R}^{m\times n}$ with $\mbox{rank}(X)=n$, \mbox{Sketch matrix} $S \in \mathbb{R}^{s\times m}$.
\ENSURE \mbox{Orthogonal factor} $Q \in \mathbb{R}^{m\times n}$, \mbox{Upper triangular factor} $R \in \mathbb{R}^{n \times n}.$ 
\STATE $K=SX,$
\STATE $[Q_{s},R_{1}]=\mbox{HouseholderQR}(K),$
\STATE $Q_{1}=XR_{1}^{-1},$
\STATE $[Q,Z]=\mbox{CholeskyQR}(Q_{1}),$
\STATE $R=ZR_{1}.$
\end{algorithmic}
\end{algorithm}%

In the following, we present the theoretical results of RCholeskyQR with the single sketch in Lemma~\ref{lemma 3.4} from \cite{Householder}.

\begin{lemma}[Rounding error analysis of RCholeskyQR]
\label{lemma 3.4}
Suppose $0 \le \epsilon<\frac{308}{317}$ and $S \in \mathbb{R}^{p\times n}$ is a $(\epsilon,p,n)$ oblivious $l_{2}$-subspace embedding. Furthermore, suppose $X \in \mathbb{R}^{m\times n}$ satisfies $\mbox{rank}(X)=n$ and $1<n \le s \le m$ where $mn\uu \le \frac{1}{12}$, $s^{\frac{3}{2}}\uu \le \frac{1}{12}$ and $\delta=\frac{383(sn^{\frac{3}{2}}+\sqrt{n}(s^{\frac{3}{2}}\sqrt{1+\epsilon}+m\norm{S}_{F})}{\sqrt{1-\epsilon}}\uu \cdot \kappa_{2}(X) \le 1$. When the above assumptions are satisfied, for $[\hat{Q},\hat{R}]=\mbox{RCholeskyQR}(X)$, we have
\begin{align}
\norm{\hat{Q}^{\top}\hat{Q}-I_{n}}_{F} &\le \frac{5445}{(25\sqrt{\frac{1-\epsilon}{1+\epsilon}}-3)^{2}}(mn\uu+n(n+1)\uu), \nonumber \\
\norm{\hat{Q}\hat{R}-X}_{F} &\le (\frac{56}{25\frac{1-\epsilon}{\sqrt{1+\epsilon}}-3\sqrt{1-\epsilon}}+\frac{1.5}{\sqrt{1-\epsilon}}\sqrt{1+\frac{5445(mn\uu+n(n+1)\uu}{(25\sqrt{\frac{1-\epsilon}{1+\epsilon}}-3)^{2}}}) \nonumber \\ &\cdot (\sqrt{1+\epsilon}\norm{X}_{2}+\frac{1-\epsilon}{12}\sigma_{n}(X)\delta)n^{2}\uu+\frac{\delta}{10}\sigma_{n}(X), \nonumber
\end{align}
with probability at least $1-p$.
\end{lemma}

\section{Numerical stability of RCLUPP}
\label{sec:RCLUPP}
In this section, we show the numerical stability of RCLUPP and some of its properties. We start with the popular randomized technique, matrix sketching \cite{rgs, tool}. For $S\in\mathbb{R}^{s\times m}$ satisfying $s\ll m$, with $0 \le \epsilon<1$, $S$ is called an $\epsilon$-embedding for an $n$-dimensional subspace $\mathcal{K}\subset\mathbb{R}^m$ when $\forall x\in\mathcal{K}$, we always have
\begin{equation} \label{eq:sketch}
(1-\epsilon)\norm{x}_2^2\leq \norm{Sx}_2^2\leq (1+\epsilon)\norm{x}_2^2.
\end{equation}
If\eqref{eq:sketch} holds for every fixed $n$-dimensional subspace of $\mathbb{R}^m$ with probability at least $1-p$, then $S$ is referred to as an $(\epsilon,p,n)$ oblivious $\ell_2$-subspace embedding. In particular, for $X \in \mathbb{R}^{m\times n}$ with full column rank, when $S \in \mathbb{R}^{s\times m}$ is an $\epsilon$-embedding for $\mbox{range}(X)$, then we have
\begin{align}
\sqrt{1-\epsilon} \cdot \norm{X}_{2} &\le \norm{SX}_{2} \le \sqrt{1+\epsilon} \cdot \norm{X}_{2}, \label{eq:s2} \\
\sqrt{1-\epsilon} \cdot \sigma_{n}(X) &\le \sigma_{n}(SX) \le \norm{SX}_{2} \le \sqrt{1+\epsilon} \cdot \norm{X}_{2}, \label{eq:sigman1} \\
\frac{\sigma_{n}(SX)}{\sqrt{1+\epsilon}} &\le \sigma_{n}(X) \le \norm{X}_{2} \le \frac{\norm{SX}_{2}}{\sqrt{1-\epsilon}}, \label{eq:sigman2}
\end{align}
with probability at least $1-p$. Common choices of the sketch matrix $S$ include the Gaussian sketch, the CountSketch and SRHT, see \cite{rgs, pmlr, Estimating, pylspack} for more details. In this work, we focus on the case of the Gaussian sketch in the numerical experiments with $S=\frac{1}{\sqrt{s}}G$, 
where $G\in\mathbb{R}^{s\times m}$ has independent standard Gaussian entries. In practice, a sketch size $s=\mathcal{O}(n)$, such as $s=n+1$ or $s=2n$, is often sufficient.

For $S$ in Algorithm~\ref{alg:RCLUPP}, if $S \in \mathbb{R}^{s\times m}$ is an $(\epsilon,p,n)$ oblivious $l_{2}$-subspace embedding in $\mathbb{R}^{m}$, we set
\begin{align}
\max (\frac{40}{3}\sqrt{mn\uu+n(n+1)\uu}+0.1, \frac{5}{39})<\frac{\sqrt{1-\epsilon}}{\sqrt{1+\epsilon}}<1, \label{eq:epsilon} \\
\sqrt{1-\epsilon}>1.02m\sqrt{n}\uu \cdot \norm{S}_{F}, \label{eq:epsilon1}
\end{align}
with $\epsilon \in (0,1)$. For the dimension of the problem, $X \in \mathbb{R}^{m\times n}$ and $S \in \mathbb{R}^{s\times m}$ with $n \le s \le m$, we let
\begin{equation}
\uu \le \min(\frac{1}{64p_{1}m\sqrt{n}}, \frac{1}{64sn\sqrt{n}}, \frac{\sqrt{1-\epsilon}}{1.02k_{0}m \cdot \norm{S}_{F}}), \label{eq:s}
\end{equation}
with $p_{1}=\max(\sqrt{n}, \norm{S}_{F})$ and $k_{0}=\frac{\norm{X}_{F}}{\norm{X}_{2}} \in (1,\sqrt{n})$. When we use the Gaussian sketch for RCLUPP with $s=\mathcal{O}(n)$ in the double precision, \eqref{eq:s} indicates that RCLUPP can deal with $X \in \mathbb{R}^{n\times n}$ with $n$ in the level of $10^{5}$ and $\uu=2^{-53}$. Rounding errors occur in each step of RCLUPP. With \eqref{eq:s}, we take
\begin{equation}
\gamma_{n}=\frac{n\uu}{1-n\uu} \le 1.02n\uu. \label{eq:gamma}
\end{equation}
With \eqref{eq:gamma}, based on [\cite{Higham}, (3.12)], [\cite{Higham}, Theorem 8.5], [\cite{Higham}, Theorem 9.3] and [\cite{Higham}, Section 18.3], we write the error matrices in the preconditioning steps of Algorithm~\ref{alg:RCLUPP} with the bounds as
\begin{align}
E_{a} &= SX-\hat{X}_{s}, \quad \abs{E_{a}} \le 1.02n\uu \cdot \abs{S}\abs{X}, \label{eq:es} \\
E_{b} &= P\hat{X}_{s}-\hat{M}\hat{N}, \quad \abs{E_{b}} \le 1.02n\uu \cdot \abs{\hat{M}}\abs{\hat{N}} \label{eq:elu} \\
E_{c} &= \hat{M}-\hat{Q}_{s}\hat{R}_{0}, \quad \abs{E_{c}} \le 1.02csn\uu \cdot \abs{\hat{M}}, \label{eq:eh} \\
E_{d} &= \hat{R}_{0}\hat{N}-\hat{R}_{1}, \quad \abs{E_{d}} \le 1.02n\uu \cdot \abs{\hat{R}_{0}}\abs{\hat{N}}, \label{eq:ey} \\
E_{e} &= \hat{Q}_{1}\hat{R}_{1}-X, \quad \abs{E_{e}} \le 1.02n\uu \cdot \abs{\hat{Q}_{1}}\abs{\hat{R}_{1}}. \label{eq:ewy} 
\end{align}
Here, we have $E_{a}, E_{b}, E_{c} \in \mathbb{R}^{s\times n}$, $E_{d} \in \mathbb{R}^{n\times n}$ and $E_{e} \in \mathbb{R}^{m\times n}$. $c$ in \eqref{eq:eh} is a positive constant. To simplify the rounding error analysis in this work, we take $c=1$. For LUP decomposition in \eqref{eq:eh}, we keep
\begin{align}
\kappa_{2}(\hat{M}) &\le \frac{1}{25.5snk_{1}\uu}, \label{eq:s3} \\
\kappa_{2}(\hat{M})\kappa_{2}(\hat{N}) &\le \frac{1}{21.4p_{2}}. \label{eq:s4}
\end{align}
For $p_{2}$ in \eqref{eq:s4}, we set $p_{2}=\max(k_{1}n\sqrt{n}\uu, \frac{k}{\sqrt{1-\epsilon}-k})$, $k=\frac{\norm{E_{a}}_{F}}{\norm{X}_{2}}$, $k_{1}=\frac{\norm{\hat{M}}_{F}}{\norm{\hat{M}}_{2}}$ and $k_{2}=\frac{\norm{\hat{N}}_{F}}{\norm{\hat{N}}_{2}}$. It is obvious that $1 \le k_{1}, k_{2} \le \sqrt{n}$. \eqref{eq:s3} and \eqref{eq:s4} reflects 
the applicability and some properties of RCLUPP, which will be discussed in Section~\ref{sec:comparison}.

\subsection{Some bounds of RCLUPP}
Based on the settings above, we focus on some bounds for the numerical stability of RCLUPP.

\subsubsection{Bounding $\norm{E_{a}}_{F}$}
Since we have defined $k=\frac{\norm{E_{a}}_{F}}{\norm{X}_{2}}$ in $p_{2}$ which occurs in \eqref{eq:s4}, we can bound $\norm{E_{a}}_{F}$ in \eqref{eq:es} as
\begin{equation}
\norm{E_{a}}_{F}=k\norm{X}_{2}. \label{eq:eaf}    
\end{equation}
In fact, when $k_{0}=\frac{\norm{X}_{F}}{\norm{X}_{2}}$, with \eqref{eq:es}, we can also bound $\norm{E_{a}}_{F}$ as
\begin{equation}
\norm{E_{a}}_{F} \le 1.02n\uu \cdot \norm{S}_{F}\norm{X}_{F} \le 1.02k_{0}n\uu \cdot \norm{S}_{F}\norm{X}_{2}. \label{eq:eaf1}
\end{equation}
\eqref{eq:eaf1} provides an deterministic bound of $\norm{E_{a}}_{F}$ with $\norm{X}_{2}$. The purpose to set $k$ is to get an alternative upper bound of $\norm{E_{a}}_{F}$ with $\norm{\hat{M}}_{2}\norm{\hat{N}}_{2}$ for the analysis afterwards.

\subsubsection{Bounding $\norm{E_{b}}_{F}$}
To bound $\norm{E_{b}}_{F}$ in \eqref{eq:elu}, when $k_{1}=\frac{\norm{\hat{M}}_{F}}{\norm{\hat{M}}_{2}}$ and $k_{2}=\frac{\norm{\hat{N}}_{F}}{\norm{\hat{N}}_{2}}$, we can bound $\norm{E_{b}}_{F}$ with \eqref{eq:elu} as
\begin{equation} \label{eq:lu}
\norm{E_{b}}_{F} \le 1.02n\uu \cdot \norm{\hat{M}}_{F}\norm{\hat{N}}_{F} \le 1.02k_{1}k_{2}n\uu \cdot \norm{\hat{M}}_{2}\norm{\hat{N}}_{2}.
\end{equation}

\subsubsection{Bounding $\norm{E_{c}}_{F}$}
To bound $\norm{E_{c}}_{F}$ in \eqref{eq:eh}, we need to focus on the singular values of the orthogonal factor $Q_{s}$ in \eqref{eq:eh} first. According to the theoretical result in [\cite{Higham}, Section 18.3], when taking the thin HouseholderQR for $\hat{M} \in \mathbb{R}^{s\times n}$ with the computed $Q$-factor $Q_{s}$, the following equality holds
\begin{equation}
\norm{\hat{Q}_{s}^{\top}\hat{Q}_{s}-I_{n}}_{F}=\mathcal{O}(sn\uu). \label{eq:orthos}
\end{equation}
$\mathcal{O}(sn\uu)$ in \eqref{eq:orthos} is a very small positive number. Therefore, to simplify the theoretical analysis, we let
\begin{equation}
\norm{\hat{Q}_{s}^{\top}\hat{Q}_{s}-I_{n}}_{2} \le \norm{\hat{Q}_{s}^{\top}\hat{Q}_{s}-I_{n}}_{F} \le 0.02. \label{eq:simplify}
\end{equation}
Using Weyl's theorem in [\cite{MatrixC}, Corollary 2.4.4] and with \eqref{eq:simplify}, we can get 
\begin{equation}
\sigma_{n}(\hat{Q}_{s}^{\top}\hat{Q}_{s}) \ge 1-0.02 \ge 0.98, \quad \norm{\hat{Q}_{s}^{\top}\hat{Q}_{s}}_{2} \le 1+0.02 \le 1.02. \label{eq:bound1}
\end{equation}
With \eqref{eq:bound1}, we bound the singular values of $Q_{s}$ as
\begin{equation}
\sigma_{n}(\hat{Q}_{s}) \ge 0.99, \quad \norm{\hat{Q}_{s}}_{2} \le 1.01. \label{eq:bh}
\end{equation}
With \eqref{eq:eh} and $c=1$, when $k_{1}=\frac{\norm{\hat{M}}_{F}}{\norm{\hat{M}}_{2}}$, we can bound $\norm{E_{c}}_{F}$ in \eqref{eq:eh} as
\begin{equation}
\norm{E_{c}}_{F} \le 1.02sn\uu \cdot \norm{\hat{M}}_{F} \le 1.02k_{1}sn\uu \cdot \norm{\hat{M}}_{2}. \label{eq:h}
\end{equation}

\subsubsection{Bounding $\norm{E_{d}}_{F}$}
For $\norm{E_{d}}_{F}$ in \eqref{eq:ey}, bounding $\norm{\hat{R}_{0}}_{F}$ is necessary. With the inequality of the matrix norms, for $\hat{Q}_{s} \in \mathbb{R}^{s\times n}$ and $\hat{R}_{0} \in \mathbb{R}^{n\times n}$, we have $\norm{\hat{Q}_{s}\hat{R}_{0}}_{F} \ge \sigma_{n}(\hat{Q}_{s})\norm{\hat{R}_{0}}_{F}$. With \eqref{eq:s}, \eqref{eq:eh}, \eqref{eq:bh} and \eqref{eq:h}, we can have
\begin{equation} \label{eq:sc}
\begin{split}
\sigma_{n}(\hat{Q}_{s})\norm{\hat{R}_{0}}_{F} &\le \norm{\hat{Q}_{s}\hat{R}_{0}}_{F} \\ &\le \norm{\hat{M}}_{F}+\norm{E_{c}}_{F} \\ &\le \norm{\hat{M}}_{F}+1.02sn\uu \cdot \norm{\hat{M}}_{F} \\ &\le 1.02\norm{\hat{M}}_{F}.
\end{split}
\end{equation}
Therefore, with \eqref{eq:bh} and \eqref{eq:sc}, when $k_{1}=\frac{\norm{\hat{M}}_{F}}{\norm{\hat{M}}_{2}}$, we can bound $\norm{\hat{R}_{0}}_{F}$ as
\begin{equation} \label{eq:r0f}
\norm{\hat{R}_{0}}_{F} \le \frac{1.02\norm{\hat{M}}_{F}}{\sigma_{n}(\hat{Q}_{s})} \le \frac{1.02\norm{\hat{M}}_{F}}{0.99} \le 1.04\norm{\hat{M}}_{F}=1.04k_{1}\norm{\hat{M}}_{2}. 
\end{equation}
Since $\norm{\hat{Q}_{s}\hat{R}_{0}}_{2} \ge \sigma_{n}(\hat{Q}_{s})\norm{\hat{R}_{0}}_{2}$ holds, similar to \eqref{eq:sc} and \eqref{eq:r0f}, we can also bound $\norm{\hat{R}_{0}}_{2}$ as
\begin{equation} \label{eq:r02}
\norm{\hat{R}_{0}}_{2} \le 1.04\norm{\hat{M}}_{2}. 
\end{equation}
With \eqref{eq:ey} and \eqref{eq:r0f}, when $k_{1}=\frac{\norm{\hat{M}}_{F}}{\norm{\hat{M}}_{2}}$ and $k_{2}=\frac{\norm{\hat{N}}_{F}}{\norm{\hat{N}}_{2}}$, we can bound $\norm{E_{d}}_{F}$ as
\begin{equation} \label{eq:edf}
\norm{E_{d}}_{F} \le 1.02n\uu \cdot \norm{\hat{R}_{0}}_{F}\norm{\hat{N}}_{F} \le 1.07k_{2}n\uu \cdot \norm{\hat{M}}_{F}\norm{\hat{N}}_{2}=1.07k_{1}k_{2}n\uu \cdot \norm{\hat{M}}_{2}\norm{\hat{N}}_{2}.
\end{equation}

\subsubsection{Bounding $\norm{\hat{R}_{1}}_{F}$}
For $\norm{\hat{R}_{1}}_{F}$ in \eqref{eq:ey}, since $k_{2} \le \sqrt{n}$, with \eqref{eq:es} and \eqref{eq:r0f}-\eqref{eq:edf}, we can bound $\norm{\hat{R}_{1}}_{F}$ as
\begin{equation}
\begin{split} \label{eq:r1f}
\norm{\hat{R}_{1}}_{F} &\le \norm{\hat{R}_{0}}_{F}\norm{\hat{N}}_{2}+\norm{E_{d}}_{F} \\ &\le 1.04\norm{\hat{M}}_{F}\norm{\hat{N}}_{2}+1.07k_{2}n\uu \cdot \norm{\hat{M}}_{F}\norm{\hat{N}}_{2} \\ &\le 1.06\norm{\hat{M}}_{F}\norm{\hat{N}}_{2}.
\end{split}
\end{equation}

\subsubsection{Bounding $\norm{\hat{R}_{1}^{-1}}_{2}$}
For $\norm{\hat{R}_{1}^{-1}}_{2}$, with the inequalities of the singular values, we can have $\norm{\hat{Q}_{s}}_{2}\sigma_{n}(\hat{R}_{0}) \ge \sigma_{n}(\hat{Q}_{s}\hat{R}_{0})$. With Weyl's theorem of the singular values, \eqref{eq:s}, \eqref{eq:eh} and \eqref{eq:h}, we can have
\begin{equation} \label{eq:ns1}
\begin{split}
\norm{\hat{Q}_{s}}_{2}\sigma_{n}(\hat{R}_{0}) &\ge \sigma_{n}(\hat{Q}_{s}\hat{R}_{0}) \\ &\ge \sigma_{n}(\hat{M})-\norm{E_{c}}_{F} \\ &\ge \sigma_{n}(\hat{M})-1.02k_{1}sn\uu \cdot \norm{\hat{M}}_{2} \\ &\ge 0.96\sigma_{n}(\hat{M}).
\end{split}
\end{equation}
Therefore, with \eqref{eq:bh} and \eqref{eq:ns1}, we can bound $\sigma_{n}(\hat{R}_{0})$ as
\begin{equation} \label{eq:ns}
\sigma_{n}(\hat{R}_{0}) \ge \frac{0.96\sigma_{n}(\hat{M})}{\norm{\hat{Q}_{s}}_{2}} \ge \frac{0.96\sigma_{n}(\hat{M})}{1.01} \ge 0.95\sigma_{n}(\hat{M}).
\end{equation}

Since $\hat{R}_{1}, \hat{R}_{0}, \hat{N} \in \mathbb{R}^{n\times n}$, we have $\sigma_{n}(\hat{R}_{0}\hat{N}) \ge \sigma_{n}(\hat{R}_{0})\sigma_{n}(\hat{N})$ based on the inequalities of the singular values. With Weyl's theorem of the singular values , \eqref{eq:ey}, \eqref{eq:s4}, \eqref{eq:edf} and \eqref{eq:ns}, when $k_{2} \le \sqrt{n}$, we can bound $\sigma_{n}(\hat{R}_{1})$ as
\begin{equation} \label{eq:ny}
\begin{split}
\sigma_{n}(\hat{R}_{1}) &\ge \sigma_{n}(\hat{R}_{0})\sigma_{n}(\hat{N})-\norm{E_{d}}_{F} \\ &\ge 0.95\sigma_{n}(\hat{M})\sigma_{n}(\hat{N})-1.07k_{1}k_{2}n\uu \cdot \norm{\hat{M}}_{2}\norm{\hat{N}}_{2} \\ &\ge 0.9\sigma_{n}(\hat{M})\sigma_{n}(\hat{N}).
\end{split}
\end{equation}
Therefore, with \eqref{eq:ny}, we can bound $\norm{\hat{R}_{1}^{-1}}_{2}$ as
\begin{equation} \label{eq:y-12}
\norm{\hat{R}_{1}^{-1}}_{2} = \frac{1}{\sigma_{n}(\hat{R}_{1})} \le \frac{1}{0.9\sigma_{n}(\hat{M})\sigma_{n}(\hat{N})} \le \frac{1.12}{\sigma_{n}(\hat{M})\sigma_{n}(\hat{N})}.
\end{equation}

\subsubsection{Bounding $\norm{\hat{N}\hat{R}_{1}^{-1}}$}
For $\norm{\hat{N}\hat{R}_{1}^{-1}}$, we focus on $\norm{\hat{R}_{0}^{-1}}$ first. With \eqref{eq:ns}, we can bound $\norm{\hat{R}_{0}^{-1}}_{2}$ as
\begin{equation} \label{eq:ns-1}
\norm{\hat{R}_{0}^{-1}}_{2} = \frac{1}{\sigma_{n}(\hat{R}_{0})} \le \frac{1.06}{\sigma_{n}(\hat{M})}.
\end{equation}

We rewrite \eqref{eq:ey} into $\hat{N}\hat{R}_{1}^{-1} = \hat{R}_{0}^{-1}\hat{R}_{1} \cdot \hat{R}_{1}^{-1}-\hat{R}_{0}^{-1} \cdot E_{d}\hat{R}_{1}^{-1} = \hat{R}_{0}^{-1}(I_{n}-E_{d}\hat{R}_{1}^{-1})$. Therefore, with \eqref{eq:s4}, \eqref{eq:edf}, \eqref{eq:y-12} and \eqref{eq:ns-1}, when $k_{2} \le \sqrt{n}$, we can bound $\norm{\hat{N}\hat{R}_{1}^{-1}}$ as
\begin{equation} \label{eq:uy-12}
\begin{split}
\norm{\hat{N}\hat{R}_{1}^{-1}}_{2} &\le \norm{\hat{R}_{0}^{-1}}_{2}(1+\norm{E_{d}}_{F}\norm{\hat{R}_{1}^{-1}}_{2}) \\ &\le \frac{1.06}{\sigma_{n}(\hat{M})} \cdot (1+1.07k_{1}k_{2}\uu \cdot \norm{\hat{M}}_{2}\norm{\hat{N}}_{2} \cdot \frac{1.12}{\sigma_{n}(\hat{M})\sigma_{n}(\hat{N})}) \\ &\le \frac{1.12}{\sigma_{n}(\hat{M})}.
\end{split}
\end{equation}

\subsubsection{Bounding $\norm{X\hat{R}_{1}^{-1}}_{2}$ and $\sigma_{n}(X\hat{R}_{1}^{-1})$}
To bound $\norm{X\hat{R}_{1}^{-1}}_{2}$, we focus on $\norm{E_{a}}_{F}$ first. With the rounding error in matrix multiplications, we can bound $\norm{E_{a}}_{F}$ in \eqref{eq:es} as $\norm{E_{a}}_{F} \le 1.02m\sqrt{n}\uu \cdot \norm{S}_{F}\norm{X}_{2}$. This is also the upper bound of $k$ in \eqref{eq:s4}. With \eqref{eq:s2} and \eqref{eq:es}, we can have
\begin{equation} \label{eq:a21}
\begin{split}
\norm{\hat{X}_{s}}_{2} \ge \norm{SX}_{2}-\norm{E_{a}}_{F} \ge \sqrt{1-\epsilon} \cdot \norm{X}_{2}-k\norm{X}_{2} = (\sqrt{1-\epsilon}-k) \cdot \norm{X}_{2},
\end{split}
\end{equation}
with probability at least $1-p$. \eqref{eq:epsilon1} guarantees $\sqrt{1-\epsilon}-k>0$. With \eqref{eq:s}, \eqref{eq:elu} and \eqref{eq:lu}, when $1 \le k_{1}, k_{2} \le \sqrt{n}$ and $s \ge n$, we can have 
\begin{equation} \label{eq:a22}
\begin{split}
\norm{\hat{X}_{s}}_{2} &= \norm{P\hat{X}_{s}}_{2} \\ &\le \norm{\hat{M}\hat{N}}_{2}+\norm{E_{b}}_{F} \\ &\le \norm{\hat{M}}_{2}\norm{\hat{N}}_{2}+1.02k_{1}k_{2}n\uu \cdot \norm{\hat{M}}_{2}\norm{\hat{N}}_{2} \\ &\le 1.02 \cdot \norm{\hat{M}}_{2}\norm{\hat{N}}_{2}.
\end{split}
\end{equation}
Since $\norm{E_{a}}_{F}=k\norm{X}_{2}$, with \eqref{eq:a21} and \eqref{eq:a22}, we can bound $\norm{E_{a}}_{F}$ with $\norm{\hat{M}}_{2}\norm{\hat{N}}_{2}$ as
\begin{equation} \label{eq:saf}
\norm{E_{a}}_{F} = k\norm{X}_{2} \le \frac{k}{\sqrt{1-\epsilon}-k} \cdot \norm{X_{s}}_{2} \le \frac{1.02k}{\sqrt{1-\epsilon}-k} \cdot \norm{\hat{M}}_{2}\norm{\hat{N}}_{2},
\end{equation}
with probability at least $1-p$. 

With \eqref{eq:es}-\eqref{eq:ey}, we can have
\begin{equation} \label{eq:psxy-1}
\begin{split}
PSX\hat{R}_{1}^{-1} &= P\hat{X}_{s}\hat{R}_{1}^{-1}-PE_{a}\hat{R}_{1}^{-1} \\ &= \hat{M}\hat{N}\hat{R}_{1}^{-1}-E_{b}\hat{R}_{1}^{-1}-PE_{a}\hat{R}_{1}^{-1} \\ &= \hat{Q}_{s}\hat{R}_{0}\hat{N}\hat{R}_{1}^{-1}-E_{c}\hat{N}\hat{R}_{1}^{-1}-E_{b}\hat{R}_{1}^{-1}-PE_{a}\hat{R}_{1}^{-1} \\ &= \hat{Q}_{s}-\hat{Q}_{s}E_{d}\hat{R}_{1}^{-1}-E_{c}\hat{N}\hat{R}_{1}^{-1}-E_{b}\hat{R}_{1}^{-1}-PE_{a}\hat{R}_{1}^{-1}.
\end{split}
\end{equation}
Therefore, with \eqref{eq:s3}, \eqref{eq:s4}, \eqref{eq:lu}, \eqref{eq:h}, \eqref{eq:edf}, \eqref{eq:y-12}, \eqref{eq:uy-12}, \eqref{eq:saf} and \eqref{eq:psxy-1}, we can have
\begin{equation} \label{eq:psxy-12}
\begin{split}
&\norm{PSX\hat{R}_{1}^{-1}-\hat{Q}_{s}}_{2} \\ \le &\norm{E_{d}}_{F}\norm{\hat{R}_{1}^{-1}}_{2}+\norm{E_{c}}_{F}\norm{\hat{N}\hat{R}_{1}^{-1}}_{2}+\norm{E_{b}}_{F}\norm{\hat{R}_{1}^{-1}}_{2}+\norm{E_{a}}_{F}\norm{\hat{R}_{1}^{-1}}_{2} \\ \le &1.07k_{1}k_{2}n\uu \cdot \norm{\hat{M}}_{2}\norm{\hat{N}}_{2} \cdot \frac{1.12}{\sigma_{n}(\hat{M})\sigma_{n}(\hat{N})}+1.02k_{1}sn\uu \cdot \norm{\hat{M}}_{2} \cdot \frac{1.12}{\sigma_{n}(\hat{M})} \\ +&1.02k_{1}k_{2}n\uu \cdot \norm{\hat{M}}_{2}\norm{\hat{N}}_{2} \cdot \frac{1.12}{\sigma_{n}(\hat{M})\sigma_{n}(\hat{N})}+\frac{1.02k}{\sqrt{1-\epsilon}-k} \cdot \norm{\hat{M}}_{2}\norm{\hat{N}}_{2} \cdot \frac{1.12}{\sigma_{n}(\hat{M})\sigma_{n}(\hat{N})} \\ \le &0.21,
\end{split}
\end{equation}
with probability at least $1-p$. With Weyl's theorem of the singular values, \eqref{eq:bh} and \eqref{eq:psxy-12}, we can have
\begin{align} 
\sigma_{n}(SX\hat{R}_{1}^{-1})=\sigma_{n}(PSX\hat{R}_{1}^{-1}) \ge \sigma_{n}(\hat{Q}_{s})-0.21 \ge 0.78, \label{eq:min} \\ \norm{SX\hat{R}_{1}^{-1}}_{2}=\norm{PSX\hat{R}_{1}^{-1}}_{2} \le \norm{\hat{Q}_{s}}_{2}+0.21 \le 1.22, \label{eq:max}
\end{align}
with probability at least $1-p$. Therefore, with \eqref{eq:sigman2}, \eqref{eq:min} and \eqref{eq:max}, we can bound $\norm{X\hat{R}_{1}^{-1}}_{2}$ and $\sigma_{n}(X\hat{R}_{1}^{-1})$ as
\begin{align} 
\sigma_{n}(X\hat{R}_{1}^{-1}) \ge \frac{\sigma_{n}(SX\hat{R}_{1}^{-1})}{\sqrt{1+\epsilon}} \ge \frac{0.78}{\sqrt{1+\epsilon}}, \label{eq:min1} \\
\norm{X\hat{R}_{1}^{-1}}_{2} \le \frac{\norm{SX\hat{R}_{1}^{-1}}_{2}}{\sqrt{1-\epsilon}} \le \frac{1.22}{\sqrt{1-\epsilon}}, \label{eq:max1}
\end{align}
with probability at least $1-p$. 

\subsection{Rounding error analysis of RCLUPP}
Based on the settings and theoretical results, we provide a rounding error analysis of RCLUPP with both orthogonality and residual.

\begin{theorem}[Rounding error analysis of RCLUPP]
\label{theorem:RCLUPP}
With \eqref{eq:epsilon}-\eqref{eq:s}, \eqref{eq:s3} and \eqref{eq:s4}, for $X \in \mathbb{R}^{m\times n}$ with $m \ge n$ and $\mbox{rank}(X)=n$, if $[\hat{Q},\hat{R}]=\mbox{RCLUPP}(X)$, we have
\begin{align}
\norm{\hat{Q}^{\top}\hat{Q}-I_{n}}_{F} &\le 6(mn\uu+n(n+1)\uu), \label{eq:ortho} \\
\norm{\hat{Q}\hat{R}-X}_{F} &\le \Delta, \label{eq:res}
\end{align}
with probability at least $1-p$. Here, we set $\Delta=\frac{4.36}{\sqrt{1-\epsilon}} \cdot hn\sqrt{n}\uu \cdot \norm{X}_{F}+\frac{3.02}{\sqrt{1-\epsilon}} \cdot hn^{2}\uu \cdot \norm{X}_{2}$ with $h=\frac{1}{\frac{0.78}{\sqrt{1+\epsilon}}-\frac{0.1}{\sqrt{1-\epsilon}}}$.
\end{theorem}
\begin{proof}
The proof of Theorem~\ref{theorem:RCLUPP} is divided into two parts, bounding $\norm{\hat{Q}^{\top}\hat{Q}-I_{n}}_{F}$ and $\norm{\hat{Q}\hat{R}-X}_{F}$. Similar to \eqref{eq:es}-\eqref{eq:ewy}, we write CholeskyQR2 after the preconditioning step with the error matrices in the beginning as shown in \eqref{eq:C1}-\eqref{eq:C8}. $\hat{R}_{2}, \hat{R}_{4} \in \mathbb{R}^{n\times n}$ in \eqref{eq:C2} and \eqref{eq:C6} are the computed $R$-factors from two iterative steps of CholeskyQR. $\hat{Q}_{2}, \hat{Q} \in \mathbb{R}^{m\times n}$ in \eqref{eq:C3} and \eqref{eq:C7} denote the computed $Q$-factor from CholeskyQR after the preconditioning step. $\hat{R}_{3}, \hat{R} \in \mathbb{R}^{m\times n}$ in \eqref{eq:C4} and \eqref{eq:C8} are the computed upper-triangular factors after iterative steps.
\begin{align}
\hat{G} &= \hat{Q}_{1}^{\top}\hat{Q}_{1}+E_{1}, \quad \abs{E_{1}} \le 1.02m\uu \cdot \abs{\hat{Q}_{1}^{\top}}\abs{\hat{Q}_{1}}, \label{eq:C1} \\
\hat{R}_{2}^{\top}\hat{R}_{2} &= \hat{G}+E_{2}, \quad \abs{E_{2}} \le 1.02(n+1)\uu \cdot \abs{\hat{R}_{2}^{\top}}\abs{\hat{R}_{2}}, \label{eq:C2}\\
\hat{Q}_{2}\hat{R}_{2} &= \hat{Q}_{1}+E_{3}, \quad \abs{E_{3}} \le 1.02n\uu \cdot \abs{\hat{Q}_{2}}\abs{\hat{R}_{2}}, \label{eq:C3} \\
\hat{R}_{3} &= \hat{R}_{2}\hat{R}_{1}+E_{4}, \quad \abs{E_{4}} \le 1.02n\uu \cdot \abs{\hat{R}_{2}}\abs{\hat{R}_{1}}, \label{eq:C4} \\
\hat{G_{1}} &= \hat{Q}_{2}^{\top}\hat{Q}_{2}+E_{5}, \quad \abs{E_{5}} \le 1.02m\uu \cdot \abs{\hat{Q}_{2}^{\top}}\abs{\hat{Q}_{2}}, \label{eq:C5} \\
\hat{R}_{4}^{\top}\hat{R}_{4} &= \hat{G_{1}}+E_{6}, \quad \abs{E_{6}} \le 1.02(n+1)\uu \cdot \abs{\hat{R}_{4}^{\top}}\abs{\hat{R}_{4}}, \label{eq:C6} \\
\hat{Q}\hat{R}_{4} &= \hat{Q}_{2}+E_{7}, \quad \abs{E_{7}} \le 1.02n\uu \cdot \abs{\hat{Q}}\abs{\hat{R}_{4}}, \label{eq:C7} \\
\hat{R} &= \hat{R}_{4}\hat{R}_{3}+E_{8}, \quad \abs{E_{8}} \le 1.02n\uu \cdot \abs{\hat{R}_{4}}\abs{\hat{R}_{3}}. \label{eq:C8}
\end{align}
Here, $Z$ in Algorithm~\ref{alg:RCLUPP} satisfies $Z=R_{4}R_{2}$ without considering the rounding errors.

\subsubsection{Bounding $\norm{\hat{Q}^{\top}\hat{Q}-I_{n}}_{F}$}
To bound $\norm{\hat{Q}^{\top}\hat{Q}-I_{n}}_{F}$ with $\hat{Q}$ in \eqref{eq:C8}, we focus on $\kappa_{2}(\hat{Q}_{1})=\frac{\norm{\hat{Q}_{1}}_{2}}{\sigma_{n}(\hat{Q}_{1})}$ first. For $\norm{\hat{Q}_{1}}_{2}$, with \eqref{eq:ewy}, \eqref{eq:s4}, \eqref{eq:r1f}, \eqref{eq:y-12}, \eqref{eq:max1}, when $k_{1}=\frac{\norm{\hat{M}}_{F}}{\norm{\hat{M}}_{2}}$, we can have
\begin{equation} \label{eq:w2}
\begin{split}
\norm{\hat{Q}_{1}}_{2} &\le \norm{X\hat{R}_{1}^{-1}}_{2}+\norm{E_{e}}_{2}\norm{\hat{R}_{1}^{-1}}_{2} \\ &\le \norm{X\hat{R}_{1}^{-1}}_{2}+1.02n\uu \cdot \norm{\hat{Q}_{1}}_{F}\norm{\hat{R}_{1}}_{F}\norm{\hat{R}_{1}^{-1}}_{2} \\ &\le \norm{X\hat{R}_{1}^{-1}}_{2}+1.02n\sqrt{n}\uu \cdot \norm{\hat{Q}_{1}}_{2}\norm{\hat{R}_{1}}_{F}\norm{\hat{R}_{1}^{-1}}_{2} \\ &\le \frac{1.22}{\sqrt{1-\epsilon}}+1.02k_{1}n\sqrt{n}\uu \cdot 1.06\norm{\hat{M}}_{2}\norm{\hat{N}}_{2} \cdot \frac{1.12}{\sigma_{n}(\hat{M})\sigma_{n}(\hat{N})} \cdot \norm{\hat{Q}_{1}}_{2} \\ &\le \frac{1.22}{\sqrt{1-\epsilon}}+0.06\norm{\hat{Q}_{1}}_{2},
\end{split}
\end{equation}
with probability at least $1-p$. Therefore, with \eqref{eq:w2}, we can bound $\norm{\hat{Q}_{1}}_{2}$ as
\begin{equation} \label{eq:w22}
\norm{\hat{Q}_{1}}_{2} \le \frac{1.3}{\sqrt{1-\epsilon}},
\end{equation}
with probability at least $1-p$. For $\sigma_{n}(\hat{Q}_{1})$, using Wely's theorem of the singular values, with \eqref{eq:ewy}, \eqref{eq:y-12}, \eqref{eq:min1}, \eqref{eq:w2} and \eqref{eq:w22}, we can bound $\sigma_{n}(\hat{Q}_{1})$ as
\begin{equation} \label{eq:nw}
\begin{split}
\sigma_{n}(\hat{Q}_{1}) &\ge \sigma_{n}(X\hat{R}_{1}^{-1})-\norm{E_{e}}_{2}\norm{\hat{R}_{1}^{-1}}_{2} \\ &\ge \sigma_{n}(X\hat{R}_{1}^{-1})-0.06\norm{\hat{Q}_{1}}_{2} \\ &\ge \frac{0.78}{\sqrt{1+\epsilon}}-0.06 \cdot \frac{1.3}{\sqrt{1-\epsilon}} \\ &\ge \frac{0.78}{\sqrt{1+\epsilon}}-\frac{0.078}{\sqrt{1-\epsilon}},
\end{split}
\end{equation}
with probability at least $1-p$.

Therefore, for $\kappa_{2}(\hat{Q}_{1})$, with \eqref{eq:epsilon}, \eqref{eq:w22} and \eqref{eq:nw}, we can have
\begin{equation} \label{eq:k2w}
\kappa_{2}(\hat{Q}_{1}) = \frac{\norm{\hat{Q}_{1}}_{2}}{\sigma_{n}(\hat{Q}_{1})} \le \frac{50}{30\sqrt{\frac{1-\epsilon}{1+\epsilon}}-3} \le \frac{1}{8\sqrt{mn\uu+n(n+1)\uu}}. 
\end{equation}
Therefore, according to \cite[Theorem 3.3]{error} and with \eqref{eq:k2w}, we can bound $\norm{\hat{Q}^{\top}\hat{Q}-I_{n}}_{F}$ as
\begin{equation} 
\norm{\hat{Q}^{\top}\hat{Q}-I_{n}}_{F} \le 6(mn\uu+n(n+1)\uu), \nonumber
\end{equation}
with probability at least $1-p$. \eqref{eq:ortho} is proved.

\subsubsection{Bounding $\norm{\hat{Q}\hat{R}-X}_{F}$}
With \eqref{eq:C1}-\eqref{eq:C8}, we can have
\begin{equation} \label{eq:qr-x}
\begin{split}
\hat{Q}\hat{R} &= (\hat{Q}_{2}+E_{7})\hat{R}_{4}^{-1}(\hat{R}_{4}\hat{R}_{3}+E_{8}) \\ &= \hat{Q}_{2}\hat{R}_{3}+E_{7}\hat{R}_{3}+\hat{Q}E_{8} \\ &= (\hat{Q}_{1}+E_{3})\hat{R}_{2}^{-1}(\hat{R}_{2}\hat{R}_{1}+E_{4})+E_{7}\hat{R}_{3}+\hat{Q}E_{8} \\ &= \hat{Q}_{1}\hat{R}_{1}+E_{3}\hat{R}_{1}+\hat{Q}_{2}E_{4}+E_{7}\hat{R}_{3}+\hat{Q}E_{8}. 
\end{split}
\end{equation}
Therefore, for $\norm{\hat{Q}\hat{R}-X}_{F}$, with \eqref{eq:qr-x}, we can have
\begin{equation} \label{eq:qr-xf}
\begin{split}
\norm{\hat{Q}\hat{R}-X}_{F} &\le \norm{\hat{Q}_{1}\hat{R}_{1}-X}_{F}+\norm{E_{3}}_{F}\norm{\hat{R}_{1}}_{2}+\norm{\hat{Q}_{2}}_{2}\norm{E_{4}}_{F} \\ &+\norm{E_{7}}_{F}\norm{\hat{R}_{3}}_{2}+\norm{\hat{Q}}_{2}\norm{E_{8}}_{F}.
\end{split}
\end{equation}

\textbf{Bounding $\norm{E_{e}}_{F}=\norm{\hat{Q}_{1}\hat{R}_{1}-X}_{F}$ with $\norm{X}_{F}$ and $\norm{\hat{R}_{1}}_{2}$} \\
For $E_{e}=\hat{Q}_{1}\hat{R}_{1}-X$ in \eqref{eq:ewy}, although we have already provided an upper bound of $\norm{E_{e}}_{F}$ with $\norm{\hat{M}}_{2}\norm{\hat{N}}_{2}$ in \eqref{eq:w2}, we try to bound $\norm{\hat{R}_{1}}_{F}$ with $\norm{X}_{F}$ for \eqref{eq:qr-xf}. With the inequalities of the singular values as used to get \eqref{eq:ewy}, \eqref{eq:sc} and \eqref{eq:w2}, we can have
\begin{equation} \label{eq:xf}
\begin{split}
\norm{\hat{Q}_{1}\hat{R}_{1}-E_{e}}_{F} &= \norm{X}_{F} \\ &\ge \norm{\hat{Q}_{1}\hat{R}_{1}}_{F}-\norm{E_{e}}_{F} \\ &\ge \sigma_{n}(\hat{Q}_{1})\norm{\hat{R}_{1}}_{F}-1.02n\sqrt{n}\uu \cdot \norm{\hat{Q}_{1}}_{2}\norm{\hat{R}_{1}}_{F}. 
\end{split}
\end{equation}
With \eqref{eq:s}, \eqref{eq:w22}, \eqref{eq:nw} and \eqref{eq:xf}, we can bound $\norm{\hat{R}_{1}}_{F}$ with $\norm{X}_{F}$ as
\begin{equation} \label{eq:yf}
\begin{split}
\norm{\hat{R}_{1}}_{F} &\le \frac{\norm{X}_{F}}{\sigma_{n}(\hat{Q}_{1})-1.02n\sqrt{n}\uu \cdot \norm{\hat{Q}_{1}}_{2}} \\ &\le \frac{\norm{X}_{F}}{\frac{0.78}{\sqrt{1+\epsilon}}-\frac{0.078}{\sqrt{1-\epsilon}}-1.02n\sqrt{n}\uu \cdot \frac{1.3}{\sqrt{1-\epsilon}}} \\ &\le \frac{\norm{X}_{F}}{\frac{0.78}{\sqrt{1+\epsilon}}-\frac{0.1}{\sqrt{1-\epsilon}}} \\ &= h\norm{X}_{F}, 
\end{split}
\end{equation}
with probability at least $1-p$. Here, $h=\frac{1}{\frac{0.78}{\sqrt{1+\epsilon}}-\frac{0.1}{\sqrt{1-\epsilon}}}$. In a similar way, we can also bound $\norm{\hat{R}_{1}}_{2}$ as
\begin{equation} \label{eq:y2a}
\norm{\hat{R}_{1}}_{2} \le h\norm{X}_{2},
\end{equation}
with probability at least $1-p$. Therefore, with \eqref{eq:ewy}, \eqref{eq:w22} and \eqref{eq:yf}, we can bound $\norm{E_{e}}_{F}=\norm{\hat{Q}_{1}\hat{R}_{1}-X}_{F}$ with $\norm{X}_{F}$ as
\begin{equation} \label{eq:wy-xf}
\begin{split}
\norm{E_{e}}_{F} &\le 1.02n\sqrt{n}\uu \cdot \norm{W}_{2}\norm{Y}_{F} \\ &\le 1.02n\sqrt{n}\uu \cdot \frac{1.3}{\sqrt{1-\epsilon}} \cdot h\norm{X}_{F} \\ &\le \frac{1.33}{\sqrt{1-\epsilon}} \cdot hn\sqrt{n}\uu \cdot \norm{X}_{F}, 
\end{split}
\end{equation}
with probability at least $1-p$.

\textbf{Bounding $\norm{\hat{Q}_{2}}_{2}$, $\norm{E_{3}}_{F}$ and $\norm{E_{4}}_{F}$} \\
For $E_{3}$ and $E_{4}$ in \eqref{eq:C3} and \eqref{eq:C4}, we need to bound $\norm{\hat{R}_{2}}_{2}$ and $\norm{\hat{Q}_{2}}_{2}$ in \eqref{eq:C2} and \eqref{eq:C3} first. Based on \cite[(3.16), (3.26)]{error}, with \eqref{eq:s} and \eqref{eq:w22}, we can bound $\norm{\hat{R}_{2}}_{2}$ and $\norm{\hat{Q}_{2}}_{2}$ as
\begin{align}  
\norm{\hat{R}_{2}}_{2} \le \sqrt{\frac{1+1.02m\uu \cdot n}{1-1.02(n+1)\uu \cdot n}} \cdot \norm{\hat{Q}_{1}}_{2} &\le \sqrt{\frac{1+\frac{1.02}{64}}{1-\frac{1.02}{64}}} \cdot \norm{\hat{Q}_{1}}_{2} \le \frac{1.33}{\sqrt{1-\epsilon}}, \label{eq:d2} \\
\norm{\hat{Q}_{2}}_{2} &\le \frac{\sqrt{69}}{8}, \label{eq:v2}
\end{align}
with probability at least $1-p$. Therefore, with \eqref{eq:C3}, \eqref{eq:C4}, \eqref{eq:yf}, \eqref{eq:d2} and \eqref{eq:v2}, we can bound $\norm{E_{3}}_{F}$ and $\norm{E_{4}}_{F}$ as
\begin{align} \label{eq:34}
\begin{split}
\norm{E_{3}}_{F} &\le 1.02n\uu \cdot \norm{\hat{Q}_{2}}_{F}\norm{\hat{R}_{2}}_{F} \\ &\le 1.02n\uu \cdot \sqrt{n}\norm{\hat{Q}_{2}}_{2} \cdot \sqrt{n}\norm{\hat{R}_{2}}_{2} \\ &\le 1.02n\uu \cdot \frac{\sqrt{69n}}{8} \cdot \frac{1.33\sqrt{n}}{\sqrt{1-\epsilon}} \\ &\le \frac{1.41}{\sqrt{1-\epsilon}} \cdot n^{2}\uu,   
\end{split}, \quad
\begin{split}
\norm{E_{4}}_{F} &\le 1.02n\uu \cdot \norm{\hat{R}_{2}}_{F}\norm{\hat{R}_{1}}_{F} \\ &\le 1.02n\uu \cdot \sqrt{n}\norm{\hat{R}_{2}}_{2} \cdot \norm{\hat{R}_{1}}_{F} \\ &\le 1.02n\uu \cdot \frac{1.33\sqrt{n}}{\sqrt{1-\epsilon}} \cdot h\norm{X}_{F} \\ &\le \frac{1.36}{\sqrt{1-\epsilon}} \cdot hn\sqrt{n}\uu \cdot \norm{X}_{F}, 
\end{split}
\end{align}
with probability at least $1-p$. 

\textbf{Bounding $\norm{\hat{R}_{3}}_{2}$, $\norm{\hat{Q}}_{2}$,  $\norm{E_{7}}_{F}$ and $\norm{E_{8}}_{F}$} \\
For $E_{7}$ and $E_{8}$ in \eqref{eq:C7} and \eqref{eq:C8}, we need to bound $\norm{\hat{R}_{3}}_{2}$, $\norm{\hat{R}_{4}}_{2}$ and $\norm{\hat{Q}}_{2}$ in \eqref{eq:C4}, \eqref{eq:C6} and \eqref{eq:C7}. For $\hat{R}_{3}$ in \eqref{eq:C4}, with \eqref{eq:s}, \eqref{eq:yf}, \eqref{eq:d2} and \eqref{eq:34}, since $\frac{\norm{X}_{F}}{\norm{X}_{2}} \le \sqrt{n}$ for $X \in \mathbb{R}^{m\times n}$ with $m \ge n$, we can bound $\norm{\hat{R}_{3}}_{2}$ as
\begin{equation} \label{eq:n2}
\begin{split} 
\norm{\hat{R}_{3}}_{2} &\le \norm{\hat{R}_{2}}_{2}\norm{\hat{R}_{1}}_{2}+\norm{E_{4}}_{F} \\ &\le \frac{1.33}{\sqrt{1-\epsilon}} \cdot h\norm{X}_{2}+\frac{1.36}{\sqrt{1-\epsilon}} \cdot hn\sqrt{n}\uu \cdot \norm{X}_{F} \\ &\le \frac{1.36}{\sqrt{1-\epsilon}} \cdot h\norm{X}_{2},
\end{split} 
\end{equation}
with probability at least $1-p$. Similarly, with
\eqref{eq:y2a}, we can bound $\norm{\hat{R}_{3}}_{F}$ as
\begin{equation} \label{eq:nf}
\begin{split} 
\norm{\hat{R}_{3}}_{F} \le \frac{1.36}{\sqrt{1-\epsilon}} \cdot h\norm{X}_{F}, 
\end{split}
\end{equation}
with probability at least $1-p$. For $\hat{R}_{4}$ in \eqref{eq:C6}, similar to \eqref{eq:d2} and with \eqref{eq:v2}, we can bound $\norm{\hat{R}_{4}}_{2}$ as
\begin{equation} \label{eq:j22}
\norm{\hat{R}_{4}}_{2} \le \sqrt{\frac{1+\frac{1}{64}}{1-\frac{1}{64}}} \cdot \norm{\hat{Q}_{2}}_{2} \le 1.06, 
\end{equation}
with probability at least $1-p$. For $\hat{Q}$ in \eqref{eq:C7}, with \eqref{eq:s} and \eqref{eq:ortho}, we can bound $\norm{\hat{Q}}_{2}$ as
\begin{equation} \label{eq:q2}
\norm{\hat{Q}}_{2} \le \sqrt{\norm{I_{n}}_{2}+6(mn\uu+n(n+1)\uu)} \le \sqrt{1+\frac{6}{64}+\frac{6}{64}} \le 1.09, 
\end{equation}
with probability at least $1-p$. Therefore, with \eqref{eq:C7}, \eqref{eq:C8} and \eqref{eq:nf}-\eqref{eq:q2}, we can bound $\norm{E_{7}}_{F}$ and $\norm{E_{8}}_{F}$ as
\begin{align} \label{eq:78}
\begin{split}
\norm{E_{7}}_{F} &\le 1.02n\uu \cdot \norm{\hat{Q}}_{F}\norm{\hat{R}_{4}}_{F} \\ &\le 1.02n\uu \cdot \sqrt{n}\norm{\hat{Q}}_{2} \cdot \sqrt{n}\norm{\hat{R}_{4}}_{2} \\ &\le 1.02n\uu \cdot 1.09\sqrt{n} \cdot 1.06\sqrt{n} \\ &\le 1.18n^{2}\uu,
\end{split} \quad
\begin{split}
\norm{E_{8}}_{F} &\le 1.02n\uu \cdot \norm{\hat{R}_{4}}_{F}\norm{\hat{R}_{3}}_{F} \\ &\le 1.02n\uu \cdot \sqrt{n}\norm{\hat{R}_{4}}_{2} \cdot \norm{\hat{R}_{3}}_{F} \\ &\le 1.02n\uu \cdot 1.06\sqrt{n} \cdot \frac{1.36}{\sqrt{1-\epsilon}} \cdot h\norm{X}_{F} \\ &\le \frac{1.48}{\sqrt{1-\epsilon}} \cdot hn\sqrt{n}\uu \cdot \norm{X}_{F}, 
\end{split}
\end{align}
with probability at least $1-p$. Therefore, we put \eqref{eq:y2a}, \eqref{eq:wy-xf}, \eqref{eq:v2}-\eqref{eq:n2}, \eqref{eq:q2} and \eqref{eq:78} into \eqref{eq:qr-xf} and we can get \eqref{eq:res} with probability at least $1-p$. \eqref{eq:res} is proved. In all, Theorem~\ref{theorem:RCLUPP} holds.
\end{proof}

\section{Comparison of the theoretical results}
\label{sec:comparison}
After the rounding error analysis of RCLUPP, we compare the theoretical results of RCLUPP and some other existing algorithms of CholeskyQR, including CholeskyQR2, SCholeskyQR3, LC2 and RCholeskyQR. We focus on the applicability and accuracy of these algorithms. 

\subsection{Comparison of applicability}
In the beginning, we focus on the applicability of CholeskyQR-type algorithms. Although RCLUPP does not seem to have a direct theoretical restriction of $\kappa_{2}(X)$, $\norm{E_{a}}_{2}$ in \eqref{eq:es} and $\norm{E_{b}}_{2}$ in \eqref{eq:elu} are always very small. Therefore, with \eqref{eq:sigman1}, we have
\begin{equation} \label{eq:xlu}
\kappa_{2}(X) \le \alpha \cdot \frac{\sqrt{1+\epsilon}}{\sqrt{1-\epsilon}} \cdot \kappa_{2}(\hat{X}_{s}) \le \alpha \beta \cdot \frac{\sqrt{1+\epsilon}}{\sqrt{1-\epsilon}} \cdot \kappa_{2}(\hat{M})\kappa_{2}(\hat{N}).
\end{equation}
Here, $\alpha$ and $\beta$ are small positive constants satisfying $\alpha \ge 1$ and $\beta \ge 1$. With \eqref{eq:s4} and \eqref{eq:xlu}, we can bound $\kappa_{2}(X)$ as
\begin{equation} \label{eq:xlu1}
\kappa_{2}(X) \le \alpha \beta \cdot \frac{\sqrt{1+\epsilon}}{\sqrt{1-\epsilon}} \cdot \kappa_{2}(\hat{M})\kappa_{2}(\hat{N}) \le \alpha \beta \cdot \frac{\sqrt{1+\epsilon}}{\sqrt{1-\epsilon}} \cdot \frac{1}{21.4p_{2}}.
\end{equation}
This is a hidden requirement of $\kappa_{2}(X)$ for RCLUPP. In the same way, according to Lemma~\ref{lemma 3.3}, we can also bound $\kappa_{2}(X)$ for LC2 as
\begin{equation} \label{eq:xlu2}
\kappa_{2}(X) \le \theta \cdot \frac{1}{64n^{2}\uu}.
\end{equation}
Here, $\theta \ge 1$ is also a small positive constant. Therefore, we show the requirements of $\kappa_{2}(X)$ in Table~\ref{tab:comparisona} based on Lemma~\ref{lemma 3.1}-Lemma~\ref{lemma 3.4}, \eqref{eq:xlu1} and \eqref{eq:xlu2} with 
\begin{equation}
K=\frac{\sqrt{1-\epsilon}}{383(sn^{\frac{3}{2}}+\sqrt{n}(s^{\frac{3}{2}}\sqrt{1+\epsilon}+m\norm{S}_{F})\uu}. \nonumber
\end{equation}

\begin{table}
\caption{Comparison of $\kappa_{2}(X)$ for $X \in \mathbb{R}^{m\times n}$}
\centering
\begin{tabular}{||c c||}
\hline
$\mbox{Algorithms}$ & $\mbox{Sufficient condition of $\kappa_{2}(X)$}$ \\
\hline
$\mbox{CholeskyQR2}$ & $\frac{1}{8\sqrt{mn\uu+n(n+1)\uu}}$ \\
\hline
$\mbox{SCholeskyQR3}$ & $\frac{1}{86p(mn\uu+n(n+1)\uu)}$ \\
\hline
$\mbox{RCholeskyQR}$ & $\frac{1}{K}$ \\
\hline
$\mbox{LC2}$ & $\mbox{$\theta \cdot \frac{1}{64n^{2}\uu}$}$ \\
\hline
$\mbox{RCLUPP}$ & $\mbox{$\alpha \beta \cdot \frac{\sqrt{1+\epsilon}}{\sqrt{1-\epsilon}} \cdot \frac{1}{21.4p_{2}}$}$ \\
\hline
\end{tabular}
\label{tab:comparisona}
\end{table}

Table~\ref{tab:comparisona} shows that RCLUPP and LC2 achieve tighter bounds of $\kappa_{2}(X)$ than CholeskyQR2 and SCholeskyQR3. The bounds of $\kappa_{2}(X)$ for RCLUPP and RCholeskyQR are all connected to the error from matrix sketching. Numerical experiments in Section~\ref{sec:numerical} confirm these findings. To compare RCLUPP with LC2, we focus on $\kappa_{2}(\hat{M})$. Although RCLUPP performs LUP decomposition on the sketched matrix $A \in \mathbb{R}^{s\times n}$ rather than on $X \in \mathbb{R}^{m\times n}$ as shown in \eqref{eq:elu}, $\kappa_{2}(X)$ and $\kappa_{2}(\hat{X}_{s})$ are generally comparable according to \eqref{eq:xlu}. The requirements on $\kappa_{2}(\hat{M})$ for the two algorithms are compared in Lemma~\ref{tab:comparisonl} with Lemma~\ref{lemma 3.3} and \eqref{eq:s3}.

\begin{table}
\caption{Comparison of $\kappa_{2}(\hat{M})$ for $X \in \mathbb{R}^{m\times n}$}
\centering
\begin{tabular}{||c c||}
\hline
$\mbox{Algorithms}$ & $\mbox{Sufficient condition of $\kappa_{2}(\hat{M})$}$ \\
\hline
$\mbox{LC2}$ & $\frac{1}{8\sqrt{mn\uu+n(n+1)\uu}}$ \\
\hline
$\mbox{RCLUPP}$ & $\frac{1}{25.5sn\sqrt{n}\uu}$ \\
\hline
\end{tabular}
\label{tab:comparisonl}
\end{table}

Table~\ref{tab:comparisonl} shows that RCLUPP admits a significantly larger tolerance on $\kappa_{2}(\hat{M})$ of order $\mathcal{O}(\uu^{-1})$ than LC2 of order $\mathcal{O}(\uu^{-1/2})$. This advantage is confirmed by the numerical examples with special structures presented in Section~\ref{sec:numerical}. Taken together, Table~\ref{tab:comparisona} and Table~\ref{tab:comparisonl} demonstrate that RCLUPP can handle more ill-conditioned matrices than CholeskyQR2, SCholeskyQR3 and LC2, underscoring the benefit of combining LUP decomposition and the thin HouseholderQR with matrix sketching. Note that the bounds in these tables are only sufficient conditions for the applicability. In the real practice, these algorithms often perform well beyond these theoretical limits. Nevertheless, a superior theoretical bound generally indicates better practical behaviors, as is evident in the numerical experiments.

\subsection{Comparison of accuracy}
In this section, we compare the accuracy of the algorithms in terms of both orthogonality $\norm{\hat{Q}^{\top}\hat{Q}-I_{n}}_{F}$ and residual $\norm{\hat{Q}\hat{R}-X}_{F}$. A comparison of the corresponding error bounds is presented in Table~\ref{tab:comparisone}, based on Lemma~\ref{lemma 3.1}-Lemma~\ref{lemma 3.4} and Theorem~\ref{theorem:RCLUPP}. Here, we have
\begin{align}
B &=6.57n^{2}\uu\norm{X}_{g}+4.81n^{2}\uu\norm{X}_{2}, \nonumber \\
O &=\frac{5445}{(25\sqrt{\frac{1-\epsilon}{1+\epsilon}}-3)^{2}}(mn\uu+n(n+1)\uu), \nonumber \\
R &=(\frac{56}{25\frac{1-\epsilon}{\sqrt{1+\epsilon}}-3\sqrt{1-\epsilon}}+\frac{1.5}{\sqrt{1-\epsilon}}\sqrt{1+\frac{5445(mn\uu+n(n+1)\uu)}{(25\sqrt{\frac{1-\epsilon}{1+\epsilon}-3)^{2}}}}) \nonumber \\ &\cdot (\sqrt{1+\epsilon}\norm{X}_{2}+\frac{1-\epsilon}{12}\sigma_{n}(X)\delta)n^{2}\uu+\frac{\delta}{10}\sigma_{n}(X), \nonumber \\
\Delta &=\frac{4.36}{\sqrt{1-\epsilon}} \cdot hn\sqrt{n}\uu \cdot \norm{X}_{F}+\frac{3.02}{\sqrt{1-\epsilon}} \cdot hn^{2}\uu \cdot \norm{X}_{2}, \nonumber
\end{align}
with $h=\frac{0.78}{\sqrt{1+\epsilon}}-\frac{0.1}{\sqrt{1-\epsilon}}$.

\begin{table}
\caption{Comparison of the accuracy for $X \in \mathbb{R}^{m\times n}$}
\centering
\begin{tabular}{||c c c||}
\hline
$\mbox{Algorithms}$ & $\norm{\hat{Q}^{\top}\hat{Q}-I_{n}}_{F}$ & $\norm{\hat{Q}\hat{R}-X}_{F}$ \\
\hline
$\mbox{CholeskyQR2}$ & $6(mn\uu+n(n+1)\uu)$ & $5n^{2}\sqrt{n}\uu\norm{X}_{2}$ \\
\hline
$\mbox{SCholeskyQR3}$ & $6(mn\uu+n(n+1)\uu)$ & $B$ \\
\hline
$\mbox{RCholeskyQR}$ & $O$ & $R$ \\
\hline
$\mbox{LC2}$ & $6.5(mn\uu+n(n+1)\uu)$ & $4.09n^{2}\uu\norm{X}_{2}$ \\
\hline
$\mbox{RCLUPP}$ & $6(mn\uu+n(n+1)\uu)$ & $\Delta$ \\
\hline
\end{tabular}
\label{tab:comparisone}
\end{table}

For the accuracy, RCholeskyQR exhibits noticeably worse orthogonality than the other algorithms, as reflected in $O$. For the residual, we derive a sharper error bound $\Delta$ for RCLUPP of order $\mathcal{O}(n\sqrt{n}\uu\norm{X}_{F}+n^{2}\uu\norm{X}_{2})$, which is a contribution of this work showing the progress in the error analysis. The other algorithms, except CholeskyQR2, share residual bounds of order $\mathcal{O}(n^{2}\uu\norm{X}_{2})$. In fact, a more refined error analysis can also sharpen the residual bound for CholeskyQR2 to the same order. Overall, RCLUPP achieves the accuracy comparable to that of existing CholeskyQR-type algorithms in both orthogonality and residual, while offering superior orthogonality to RCholeskyQR. These theoretical findings are consistent with the numerical results in Section~\ref{sec:numerical}.

\subsection{Discussions of \eqref{eq:s4}}
\eqref{eq:s4} provides a sufficient condition for the numerical stability of RCLUPP and, through \eqref{eq:xlu1}, characterizes its applicability to $X \in \mathbb{R}^{m\times n}$. However, this condition may be conservative relative to the observed values of $\kappa_{2}(\hat{M})\kappa_{2}(\hat{N})$ in the real implementation. Compared with the condition for LC2 in \eqref{eq:xlu2}, we use $k_{1}$ and $k_{2}$ in \eqref{eq:s4}, which are typically much smaller than $\sqrt{n}$, leading to a sharper bound in \eqref{eq:s4}. Furthermore, according to the rule of thumb in \cite{Higham}, the actual rounding errors in \eqref{eq:es}-\eqref{eq:ewy} are often smaller than their most pessimistic bounds by a square root. Thus, the upper bound for $\kappa_{2}(\hat{M})\kappa_{2}(\hat{N})$ may be taken as $\frac{\sqrt{n}}{21.4p_{2}}$, which is expected to reflect the  practical behavior better. Numerical experiments concerning this sufficient condition are presented in Section~\ref{sec:numerical}.

\subsection{Some extensions-RCLUPPr}
A natural variant of RCLUPP, RCLUPPr, integrates matrix sketching with LUP decomposition and CholeskyQR in an alternative manner. In RCLUPPr, LUP decomposition is performed directly on the input matrix $X \in \mathbb{R}^{m\times n}$, after which matrix sketching is applied to the resulting $L$-factor as shown in Algorithm~\ref{alg:RCLUPPr}. A rounding error analysis similar to that in this work can be developed. Compared with RCLUPP, RCLUPPr avoids the influence of the error in matrix sketching and admits a substantially milder condition $\kappa_{2}(\hat{M})\kappa_{2}(\hat{N})=\mathcal{O}(\frac{1}{n^{2}\uu})$, which is considerably optimal than the corresponding bounds for RCholeskyQR and RCLUPP. Consequently, RCLUPPr can handle more ill-conditioned matrices, as confirmed by the numerical experiments in Section~\ref{sec:numerical}. However, performing LUP decomposition on $X$ directly incurs higher computational cost than on the compressed matrix. Generally speaking, RCLUPP achieves a superior balance among the accuracy, applicability, and efficiency for most of the problems, while RCLUPPr is recommended for the very ill-conditioned cases.

\begin{algorithm}
\caption{$[Q,R]=\mbox{RCLUPPr}(X)$}
\label{alg:RCLUPPr}
\begin{algorithmic}[1]
\REQUIRE $X \in \mathbb{R}^{m\times n}$ with $\mbox{rank}(X)=n$, \mbox{Sketch matrix} $S \in \mathbb{R}^{s\times m}$.
\ENSURE \mbox{Orthogonal factor} $Q \in \mathbb{R}^{m\times n}$, \mbox{Upper triangular factor} $R \in \mathbb{R}^{n \times n}.$
\STATE $[M,N,P]=\mbox{LU}(X),$
\STATE $M_{s}=SM,$
\STATE $[Q_{s},R_{0}]=\mbox{HouseholderQR}(L_{s}),$
\STATE $R_{1}=R_{0}N,$
\STATE $Q_{1}=XR_{1}^{-1},$
\STATE $[Q,Z]=\mbox{CholeskyQR2}(Q_{1}),$
\STATE $R=ZR_{1}.$
\end{algorithmic}
\end{algorithm}%

\section{Numerical experiments}
\label{sec:numerical}
In this section, we present numerical experiments conducted in MATLAB to evaluate the applicability, accuracy and efficiency of RCLUPP and other CholeskyQR-type algorithms. The superiority of the sharper bound of residual for RCLUPP is also demonstrated numerically. Matrices for the experiments include examples in the industry as well as those from the existing literature.

\subsection{Tests of the applicability of the algorithms}
In this section, we evaluate the applicability of the  algorithms through numerical experiments in MATLAB. We compare RCLUPP and RCLUPPr with CholeskyQR2, SCholeskyQR3, LC2 and RCholeskyQR, with particular emphasis on the ill-conditioned $X$ with large $\kappa_{2}(X)$. For the randomized algorithms, RCLUPP, RCLUPPr and RCholeskyQR, we set $\epsilon=0.5$ and $p=0.6$. \eqref{eq:epsilon} holds whenever $mn \le 10^{12}$. For SCholeskyQR3, we use the improved $s=11(mn\uu+(n+1)\uu)\norm{X}_{g}^{2}$ as shown in Lemma~\ref{lemma 3.2}. Hundreds of test cases are performed for each algorithm and the average results are reported. We assess the accuracy via the orthogonality $\norm{\hat{Q}^{\top}\hat{Q}-I_{n}}_{F}$ and residual $\norm{\hat{Q}\hat{R}-X}_{F}$.

\subsubsection{Matrix based on SVD}
\begin{figure}
\centering
\begin{minipage}{0.48\textwidth}
\centering
\includegraphics[width=\textwidth]{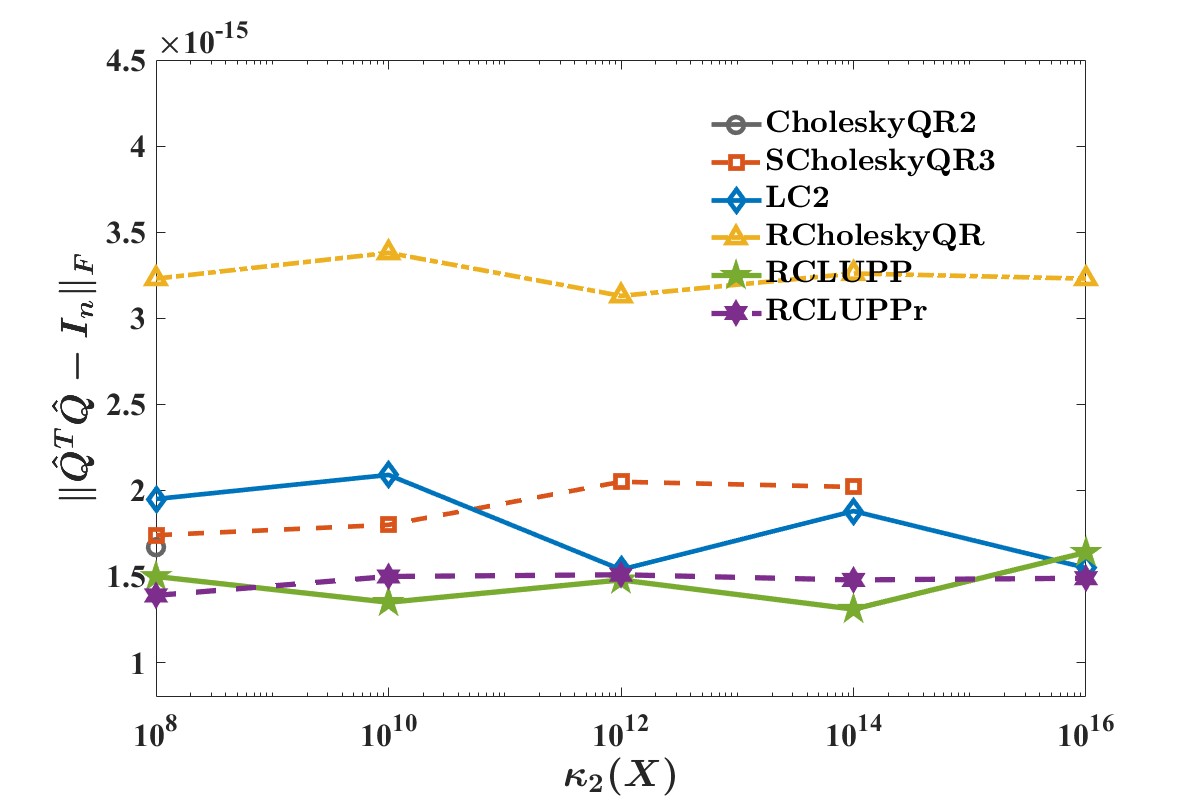}
\caption{Comparison of orthogonality on the SVD-based matrix}
\label{fig:Ok}
\end{minipage}
\hfill
\begin{minipage}{0.48\textwidth}
\centering
\includegraphics[width=\textwidth]{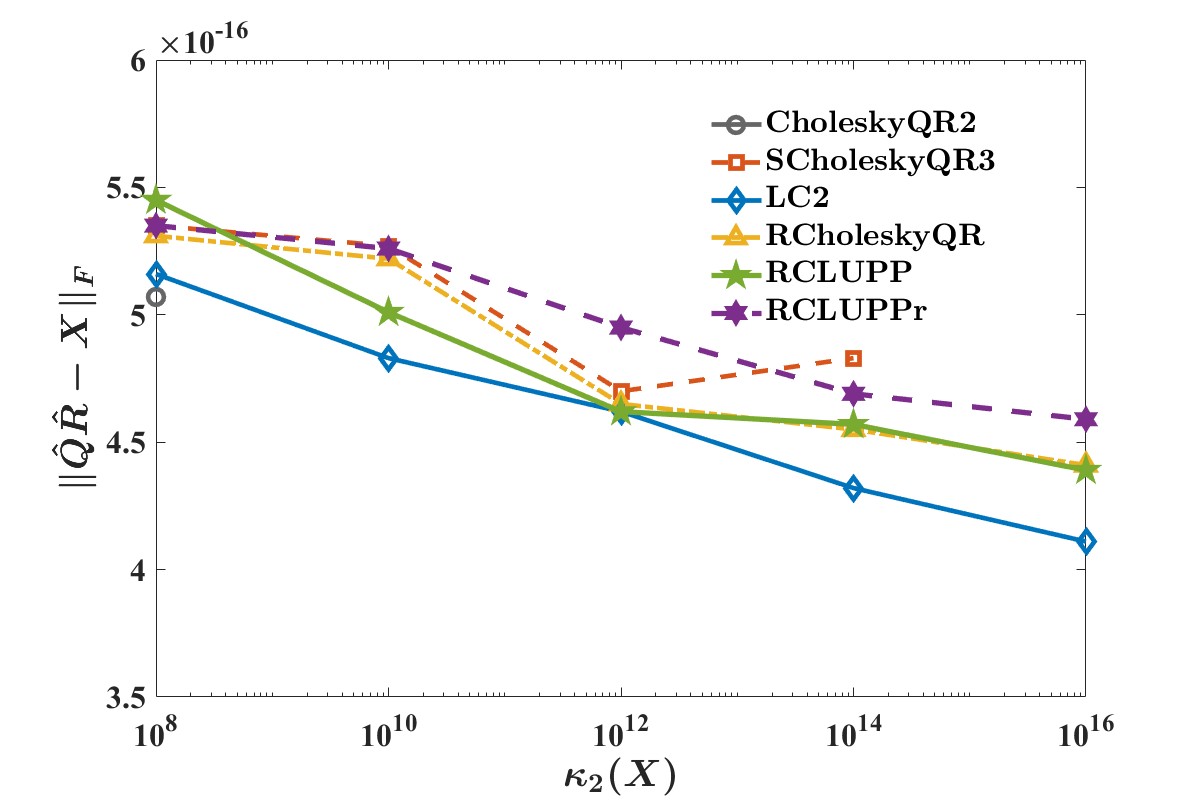}
\caption{Comparison of residual on the SVD-based matrix}
\label{fig:Rk}
\end{minipage}
\end{figure}

In this section, we focus on the comparison of applicability on the matrix based on SVD. We construct $X \in \mathbb{R}^{2000\times 50}$  using SVD. We take $s=100$ for RCLUPP, RCLUPPr and RCholeskyQR. We set
\begin{equation}
X=F \Sigma H^{\top}, \nonumber
\end{equation}
where $F \in \mathbb{R}^{2000\times 50}$ and $H \in \mathbb{R}^{50\times 50}$ are orthogonal matrices and 
\begin{equation}
\Sigma = {\rm diag}(1, \sigma^{\frac{1}{49}}, \cdots, \sigma^{\frac{48}{49}}, \sigma) \in \mathbb{R}^{50\times 50}. \nonumber
\end{equation}
Here, $0<\sigma<1$ is a positive constant. Therefore, we can have $\sigma_{1}(X)=\norm{X}_{2}=1$ and $\kappa_{2}(X)=\frac{1}{\sigma}$. We vary $\sigma$ from $10^{-10}$, $10^{-12}$, $10^{-14}$ to $10^{-16}$ to change $\kappa_{2}(X)$. Numerical results of matrix based on SVD are shown in Figure~\ref{fig:Ok} and Figure~\ref{fig:Rk}. 

Figure~\ref{fig:Ok} and Figure~\ref{fig:Rk} show that CholeskyQR2 and SCholeskyQR3 have significantly poorer applicability than LC2, RCholeskyQR, RCLUPP and RCLUPPr, consistent with the theoretical comparison in Table~\ref{tab:comparisona}. With $\kappa_{2}(X)$ approaching $\uu^{-1}$, RCholeskyQR exhibits inferior orthogonality compared with other algorithms. Consequently, RCLUPP and RCLUPPr achieve better applicability for the ill-conditioned matrices than CholeskyQR2 and SCholeskyQR3, while offering numerical stability comparable to LC2 and superior to RCholeskyQR.

\subsubsection{Lower-triangular matrix}
\begin{figure}
\centering
\begin{minipage}{0.48\textwidth}
\centering
\includegraphics[width=\textwidth]{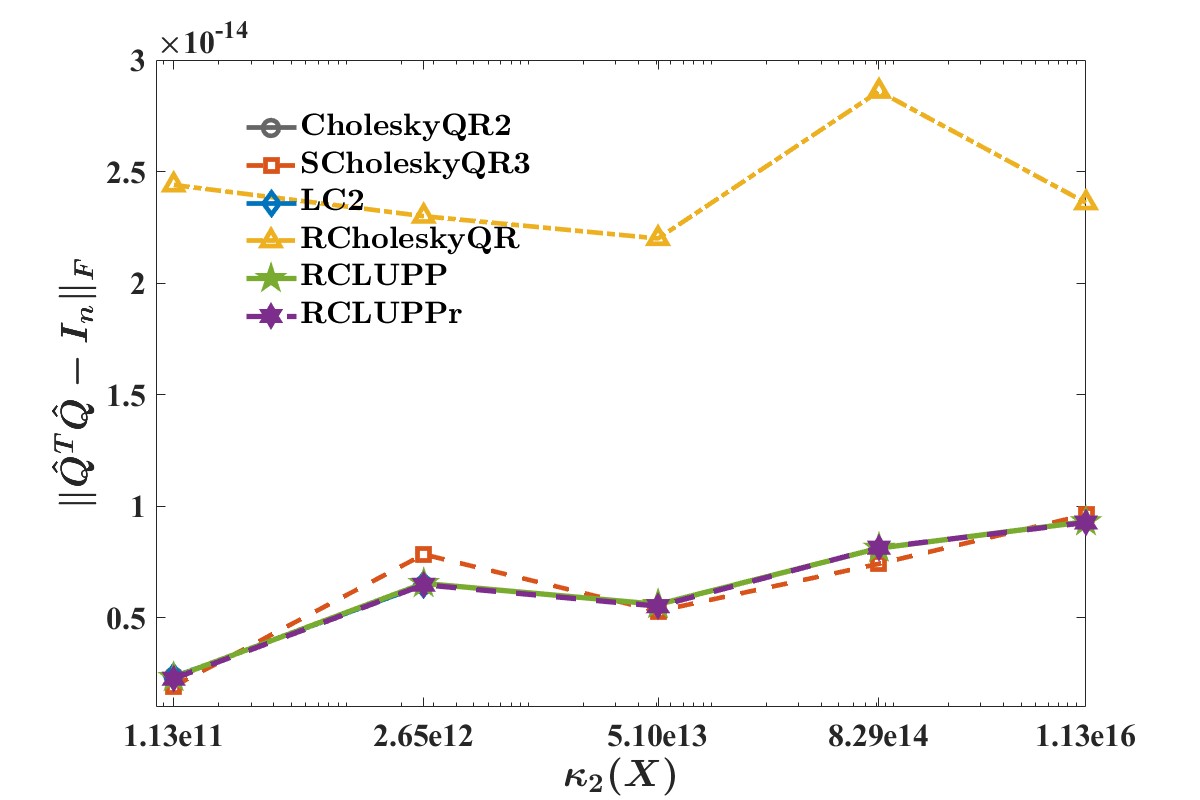}
\caption{Comparison of orthogonality on the lower-triangular matrix}
\label{fig:OK1}
\end{minipage}
\hfill
\begin{minipage}{0.48\textwidth}
\centering
\includegraphics[width=\textwidth]{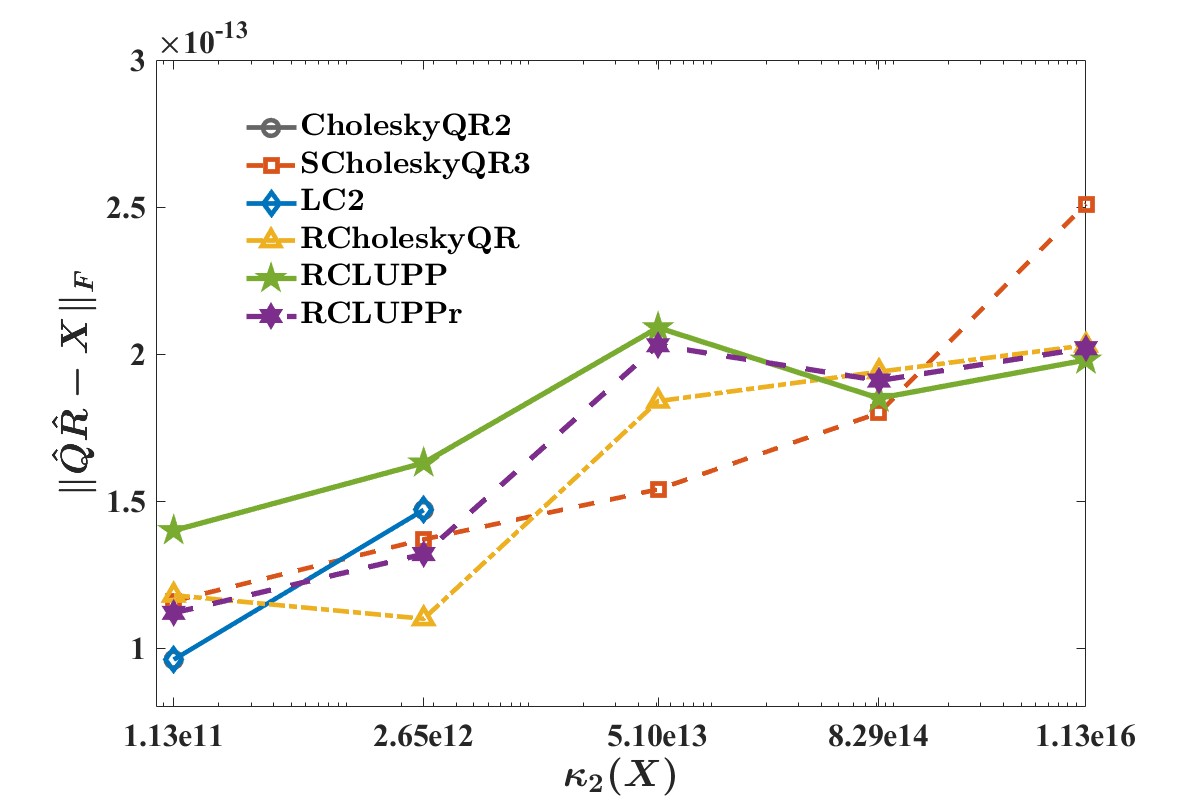}
\caption{Comparison of residual on the lower-triangular matrix}
\label{fig:RK1}
\end{minipage}
\end{figure}

Here, we consider a type of lower-triangular matrices, which occurs primarily in the numerical PDEs, financial mathematics and solving linear systems, see \cite{Theory, Geometric} and their references. We construct $X \in \mathbb{R}^{20000\times 50}$ by stacking $400$ $F \in \mathbb{R}^{50\times 50}$ from top to bottom in a block version, where $F \in \mathbb{R}^{50\times 50}$ is defined as
$F_{ij}=
\begin{cases} 
100, & \text{if } i=j, i=1,2,\cdots,50 \\
a, & \text{if } j<i, i=1,2,\cdots,50 \\
0, & others
\end{cases}.$
Here, $a$ is a negative constant. The $L$-factor after LUP decomposition of such an $X$ is ill-conditioned. In the numerical experiments, we vary $a$ from $-0.6$, $-0.7$, $-0.8$, $-0.9$ to $-1$ to alter $\kappa_{2}(X)$. As $\abs{a}$ decreases, $\kappa_{2}(X)$ increases. We take $s=100$ and the results of the numerical experiments are summarized in Figure~\ref{fig:OK1} and Figure~\ref{fig:RK1}. 

When the $L$-factor after LUP decomposition of such an $X$ is ill-conditioned, both RCLUPP and RCLUPPr show significant advantages over LC2 in dealing with this type of cases as reflected in Figure~\ref{fig:OK1} and Figure~\ref{fig:RK1},  which corresponds to Table~\ref{tab:comparisonl}. RCLUPP and RCLUPPr have the same level of accuracy in both orthogonality and residual compared with that of LC2, demonstrating that our algorithms overcome the primary problem of LC2 and are better algorithms to some extent.

\subsubsection{Arrowhead matrix}
\begin{figure}
\centering
\begin{minipage}{0.48\textwidth}
\centering
\includegraphics[width=\textwidth]{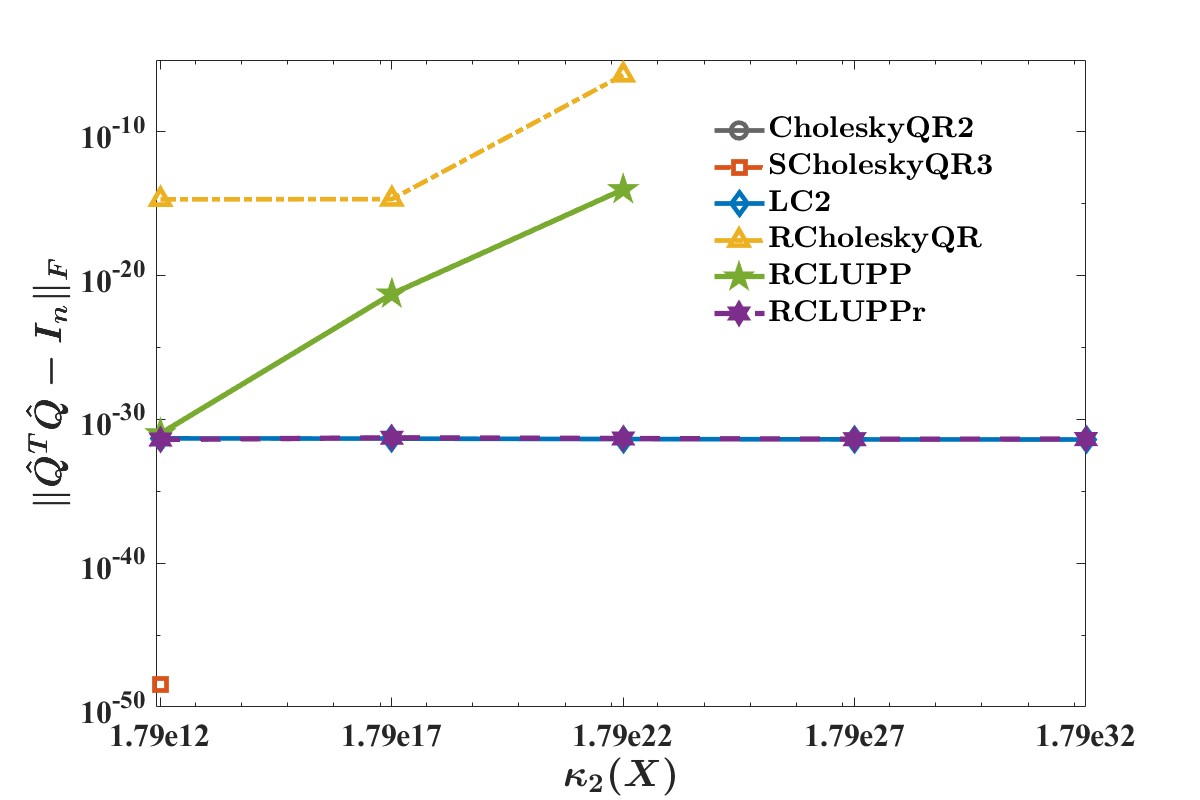}
\caption{Comparison of orthogonality on the arrowhead matrix}
\label{fig:Ok2}
\end{minipage}
\hfill
\begin{minipage}{0.48\textwidth}
\centering
\includegraphics[width=\textwidth]{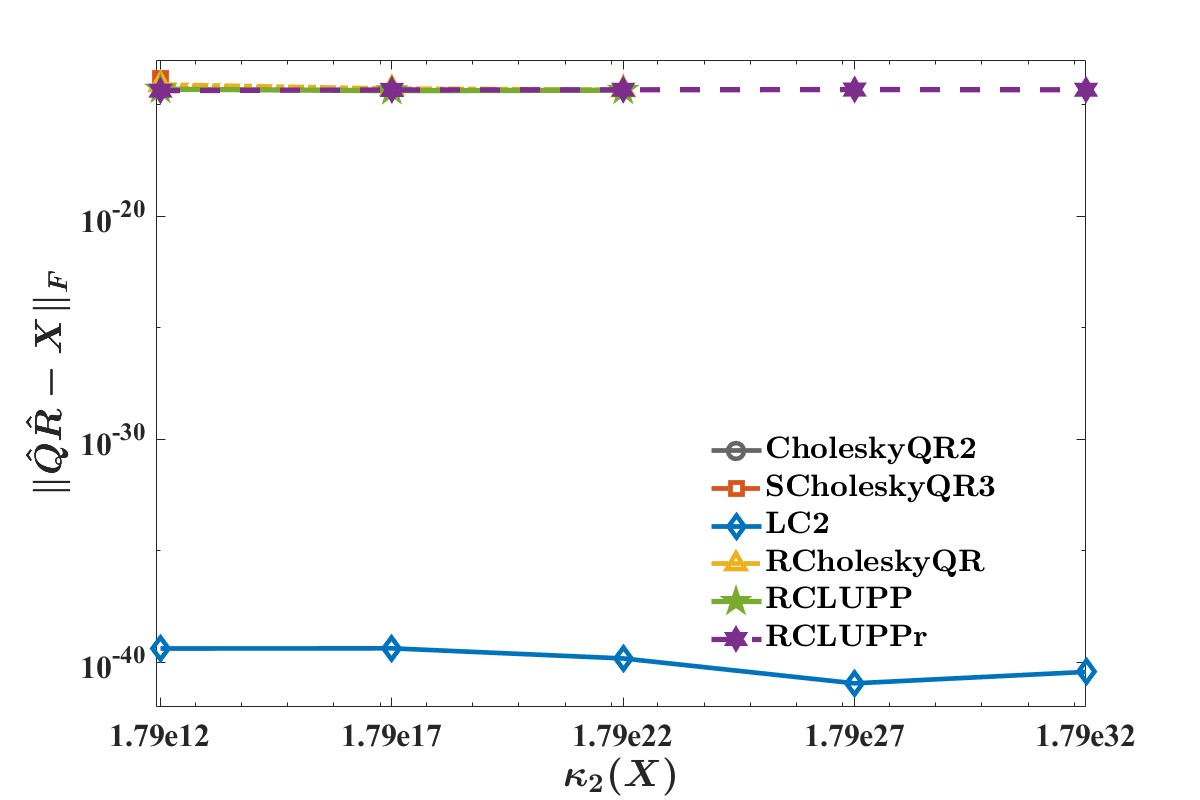}
\caption{Comparison of residual on the arrowhead matrix}
\label{fig:Rk2}
\end{minipage}
\end{figure}

In this part, we focus on the arrowhead matrix, which is widely used in graph theory, control theory, and certain eigenvalue problems \cite{Li, Eigen}. Many arrowhead matrices are very ill-conditioned. We define $e_{a}=(0,1,1, \cdots, 1,1)^{\top} \in \mathbb{R}^{50}$ and $e_{b}=(1,0,0, \cdots, 0,0)^{\top} \in \mathbb{R}^{20000}$, together with a diagonal matrix $E={\rm diag}(1, 1, \cdots, 1, \beta) \in \mathbb{R}^{50\times 50}$. Here, $\beta$ is a positive constant. Moreover, a large matrix $\mathbb{O}_{19950\times 50}$ is formed with all the elements $0$. Therefore, a matrix $P_{s} \in \mathbb{R}^{20000\times 50}$ is constructed as
\begin{equation}
P_{s}=
\begin{pmatrix}
E \\
\mathbb{O}_{19950\times 50} \nonumber
\end{pmatrix}.
\end{equation} 
Therefore, we build $X \in \mathbb{R}^{20000\times 50}$ in the form of
\begin{equation}
X=-5e_{b} \cdot e_{a}^{\top}+P_{s}. \nonumber
\end{equation}
In this group of numerical experiments, we vary $\beta$ from $10^{-10}$, $10^{-15}$, $10^{-20}$, $10^{-25}$ to $10^{-30}$ to change $\kappa_{2}(X)$. With $\beta$ decreasing, $\kappa_{2}(X)$ gets increasing. We take $s=100$ and list the numerical results in Figure~\ref{fig:Ok2} and Figure~\ref{fig:Rk2}. 

In Figure~\ref{fig:Ok2} and Figure~\ref{fig:Rk2}, RCLUPP exhibits better orthogonality with the structure of CholeskyQR2 after the preconditioning step compared with RCholeskyQR, which aligns with the comparison of the theoretical results in Table~\ref{tab:comparisone}. For the very ill-conditioned $X$, the orthogonality of RCholeskyQR decreases rapidly and the algorithm will encounter numerical breakdown for some severely ill-conditioned cases, while RCLUPP has better performance. With $\kappa_{2}(X)$ increasing, RCLUPPr shows its outstanding applicability and accuracy among all the algorithms. Based on Figure~\ref{fig:Ok}-Figure~\ref{fig:Rk2}, we conclude that RCLUPPr is more suitable for the very ill-conditioned cases. More experiments can be taken by the readers for practice.

\subsection{Tests of the accuracy of the algorithms}
\begin{figure}
\centering
\begin{minipage}{0.48\textwidth}
\centering
\includegraphics[width=\textwidth]{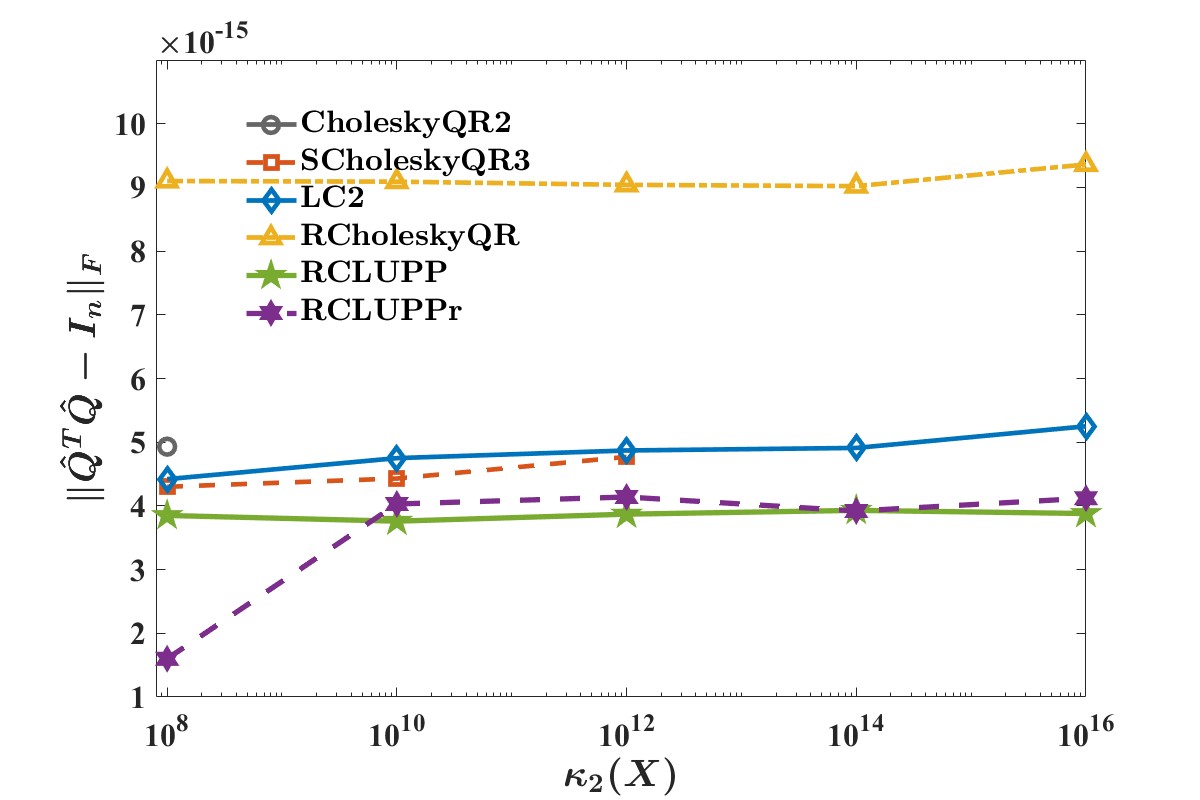}
\caption{Comparison of orthogonality with different $\kappa_{2}(X)$}
\label{fig:ko}
\end{minipage}
\hfill
\begin{minipage}{0.48\textwidth}
\centering
\includegraphics[width=\textwidth]{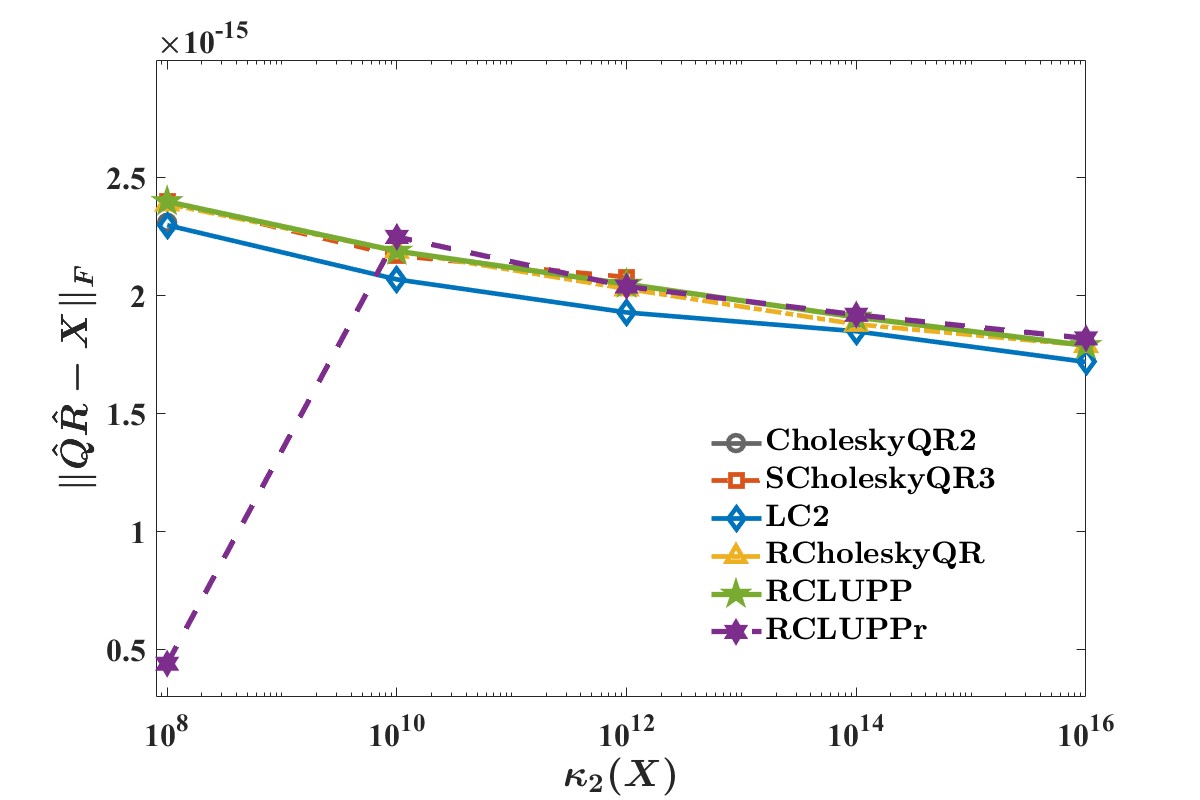}
\caption{RComparison of residual with different $\kappa_{2}(X)$}
\label{fig:kr}
\end{minipage}
\end{figure}

\begin{figure}
\centering
\begin{minipage}{0.48\textwidth}
\centering
\includegraphics[width=\textwidth]{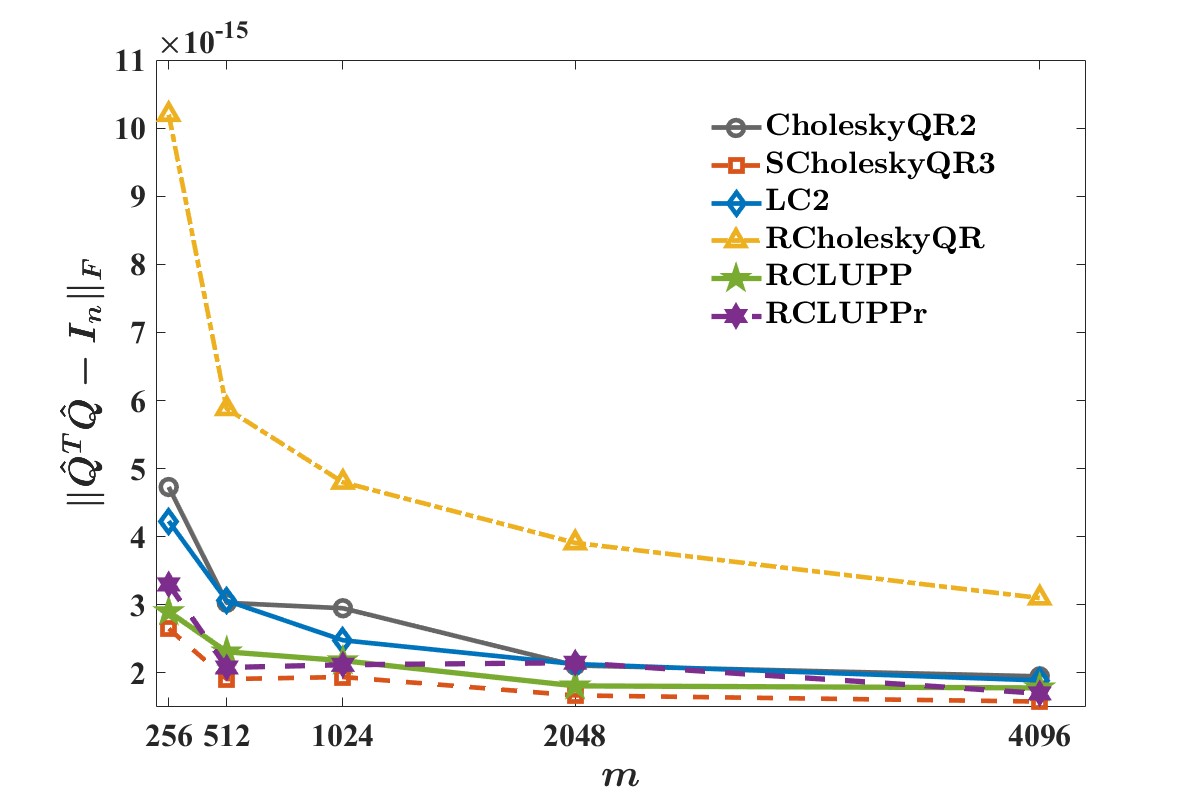}
\caption{Comparison of orthogonality with different $m$}
\label{fig:mo}
\end{minipage}
\hfill
\begin{minipage}{0.48\textwidth}
\centering
\includegraphics[width=\textwidth]{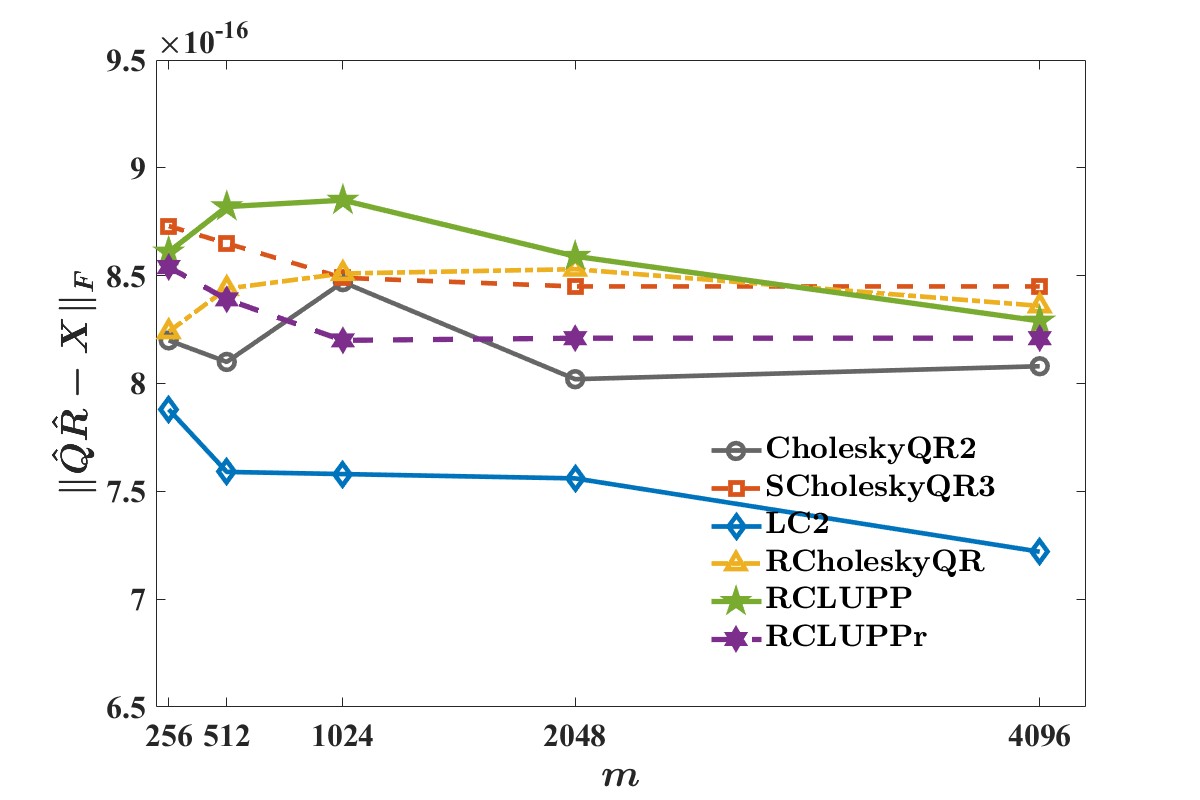}
\caption{Comparison of residual with different $m$}
\label{fig:mr}
\end{minipage}
\end{figure}

\begin{figure}
\centering
\begin{minipage}{0.48\textwidth}
\centering
\includegraphics[width=\textwidth]{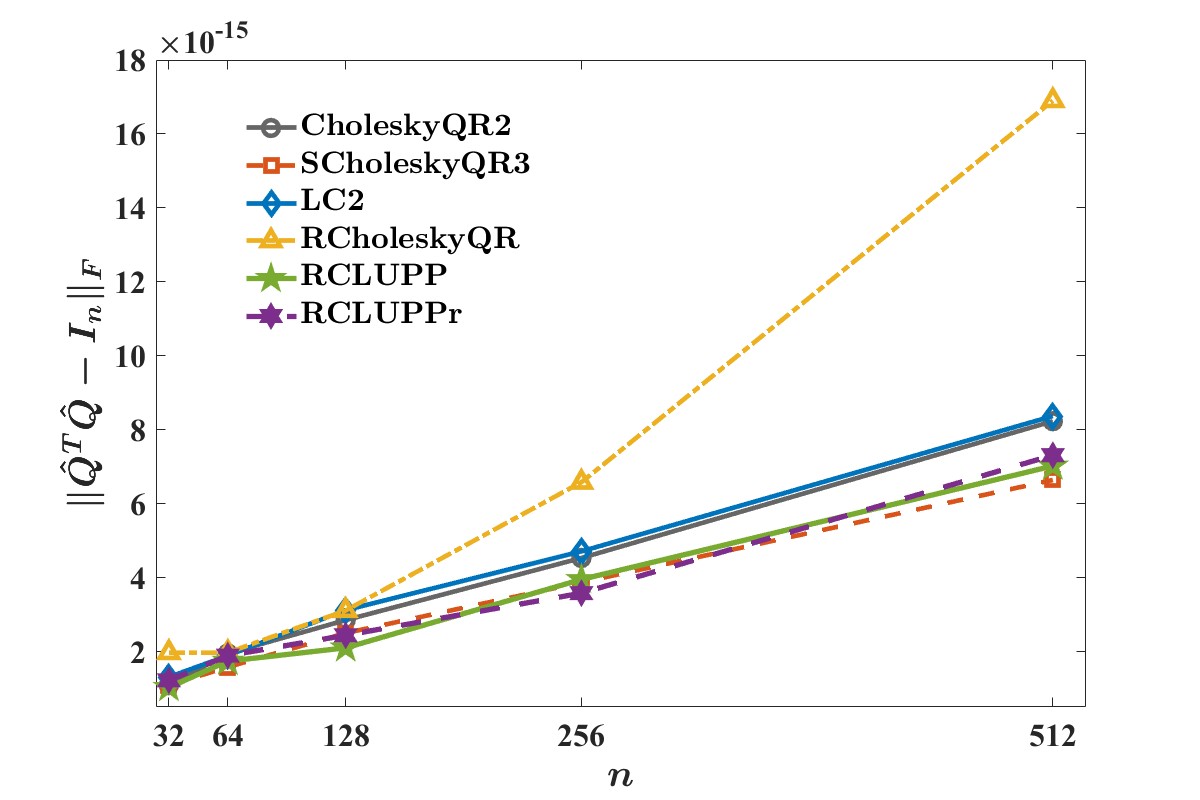}
\caption{Comparison of orthogonality between different $n$}
\label{fig:no}
\end{minipage}
\hfill
\begin{minipage}{0.48\textwidth}
\centering
\includegraphics[width=\textwidth]{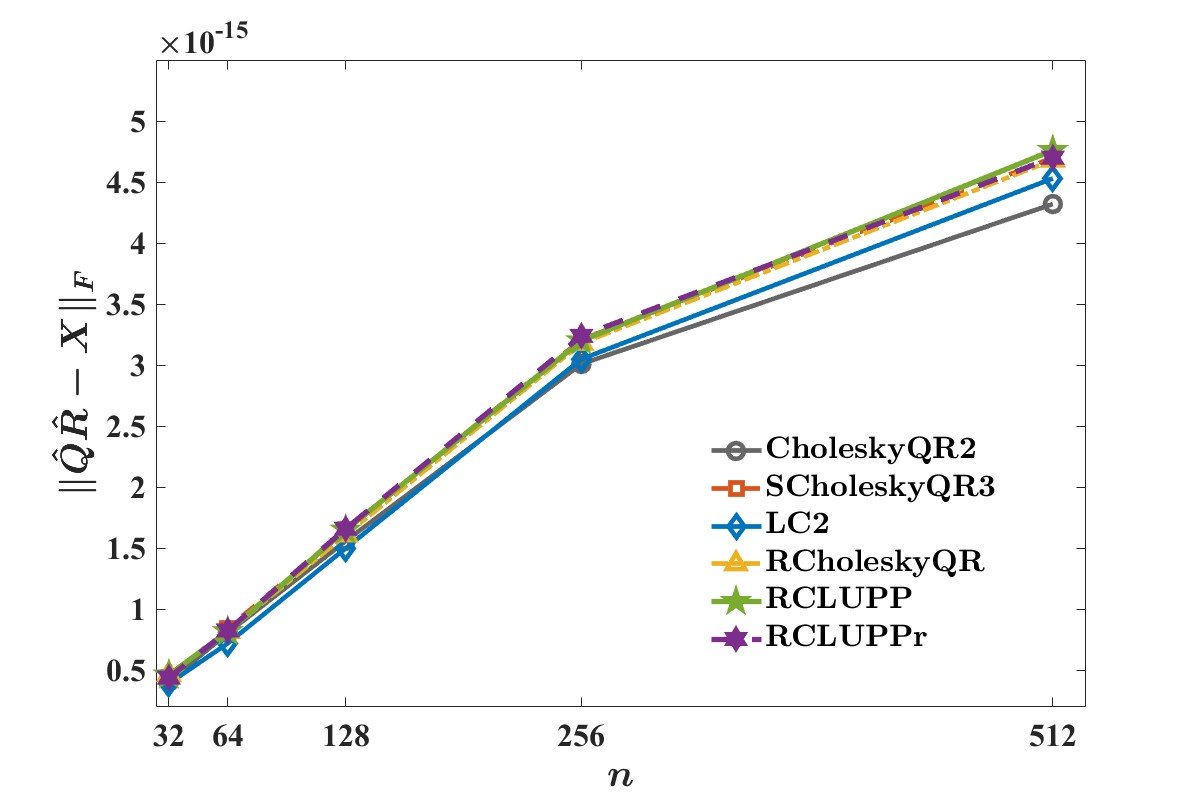}
\caption{Comparison of residual between different $n$}
\label{fig:nr}
\end{minipage}
\end{figure}

\begin{figure}
\centering
\begin{minipage}{0.48\textwidth}
\centering
\includegraphics[width=\textwidth]{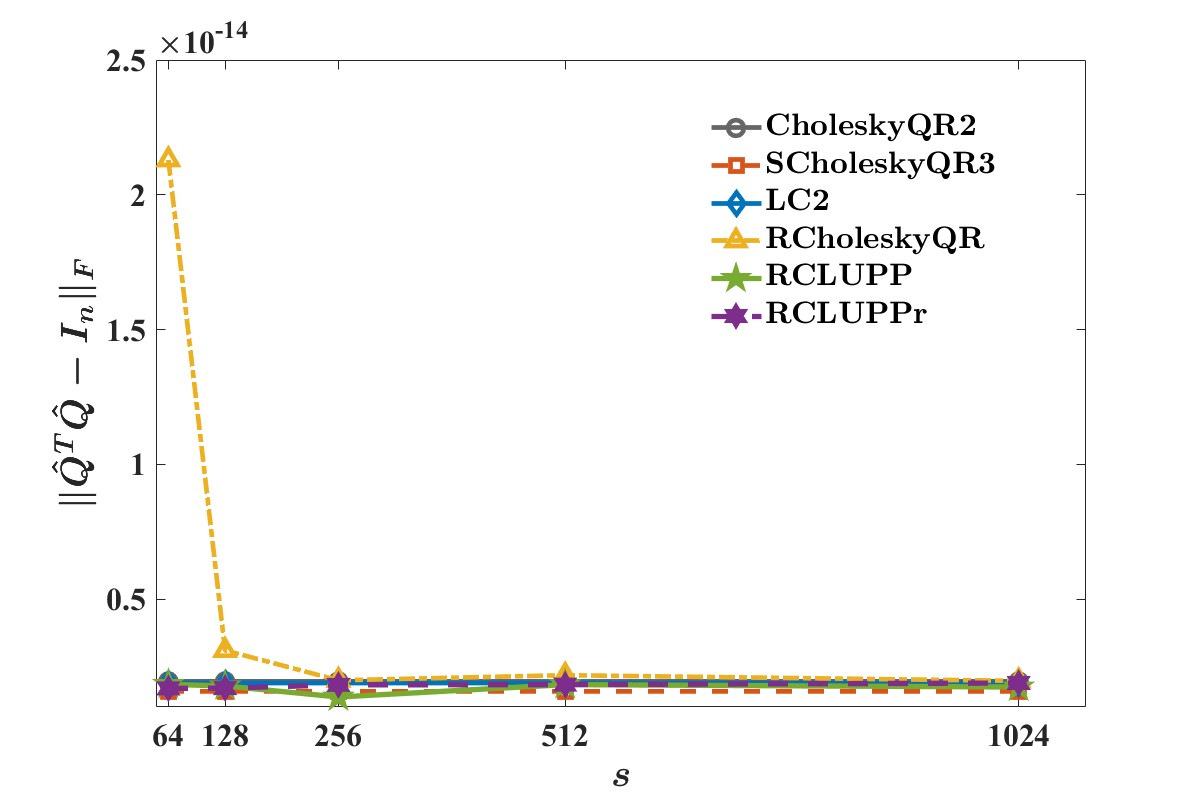}
\caption{Comparison of orthogonality between different $s$}
\label{fig:so}
\end{minipage}
\hfill
\begin{minipage}{0.48\textwidth}
\centering
\includegraphics[width=\textwidth]{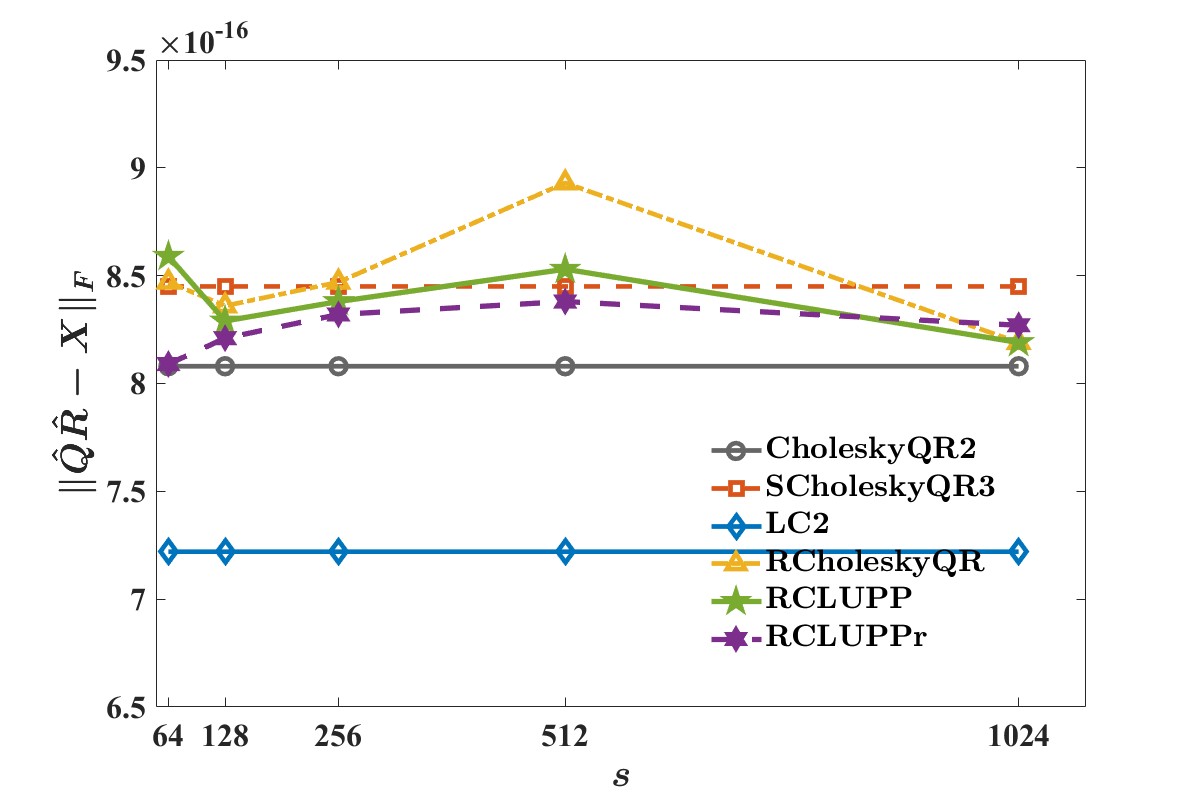}
\caption{Comparison of residual between different $s$}
\label{fig:sr}
\end{minipage}
\end{figure}

In this section, we examine the accuracy and numerical stability of the algorithms. For RCLUPP, RCLUPPr and RCholeskyQR, we set $\epsilon = 0.5$ and $p = 0.6$. We investigate the effects of $\kappa_{2}(X)$, $m$, $n$ and $s$ on the performance of the algorithms using the SVD-based matrices. To study the influence of $\kappa_{2}(X)$, we fix $m=4096$, $n=256$ and $s=512$ for the randomized algorithms. When varying $m$, we set $\kappa_{2}(X)=10^{4}$, $n=64$ and $s=128$. For different $n$, we set $\kappa_{2}(X)=10^{4}$, $m=4096$ and $s=1024$. When varying $s$, we set $\kappa_{2}(X)=10^{4}$, $m=4096$ and $n=64$. The orthogonality $\norm{\hat{Q}^{\top}\hat{Q}-I_{n}}_{F}$ and the residual $\norm{\hat{Q}\hat{R}-X}_{F}$ are reported in Figure~\ref{fig:ko}-Figure~\ref{fig:sr}.

Figure~\ref{fig:ko} and Figure~\ref{fig:kr} confirm the limited applicability of CholeskyQR2 and SCholeskyQR3, consistent with the results in Table~\ref{tab:comparisona}. CholeskyQR2 fails for the dense $X$ with $\kappa_{2}(X) \ge \uu^{-\frac{1}{2}}$, while SCholeskyQR3 becomes inapplicable as $\kappa_{2}(X)$ approaches $\uu^{-1}$. In Figure~\ref{fig:ko}-Figure~\ref{fig:sr}, RCLUPP and RCLUPPr demonstrate superior numerical stability, especially orthogonality, compared with RCholeskyQR. This aligns with that in Table~\ref{tab:comparisone}. These advantages become more pronounced for larger $\kappa_{2}(X)$. Figure~\ref{fig:so} and Figure~\ref{fig:sr} further show that the choice of $s$ has negligible impact on the accuracy of RCLUPP and RCLUPPr. Generally speaking, both algorithms maintain good accuracy and numerical stability with different $\kappa_{2}(X)$, $m$, $n$, and $s$, achieving performances comparable to other CholeskyQR-type algorithms.

\subsection{Tests of the sufficient condition}
\begin{figure}
\centering
\begin{minipage}{0.44\textwidth}
\centering
\includegraphics[width=\textwidth]{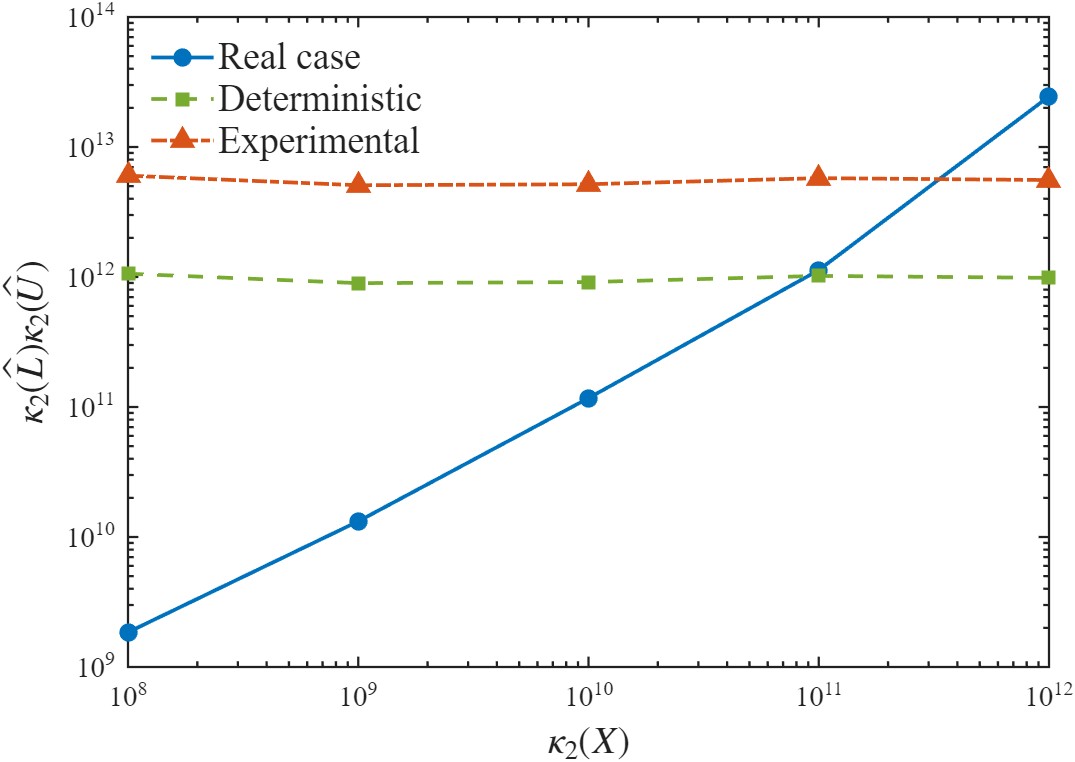}
\caption{Comparison of $\kappa_{2}(\hat{L})\kappa_{2}(\hat{U})$ with different $\kappa_{2}(X)$}
\label{fig:luk}
\end{minipage}
\hfill
\begin{minipage}{0.44\textwidth}
\centering
\includegraphics[width=\textwidth]{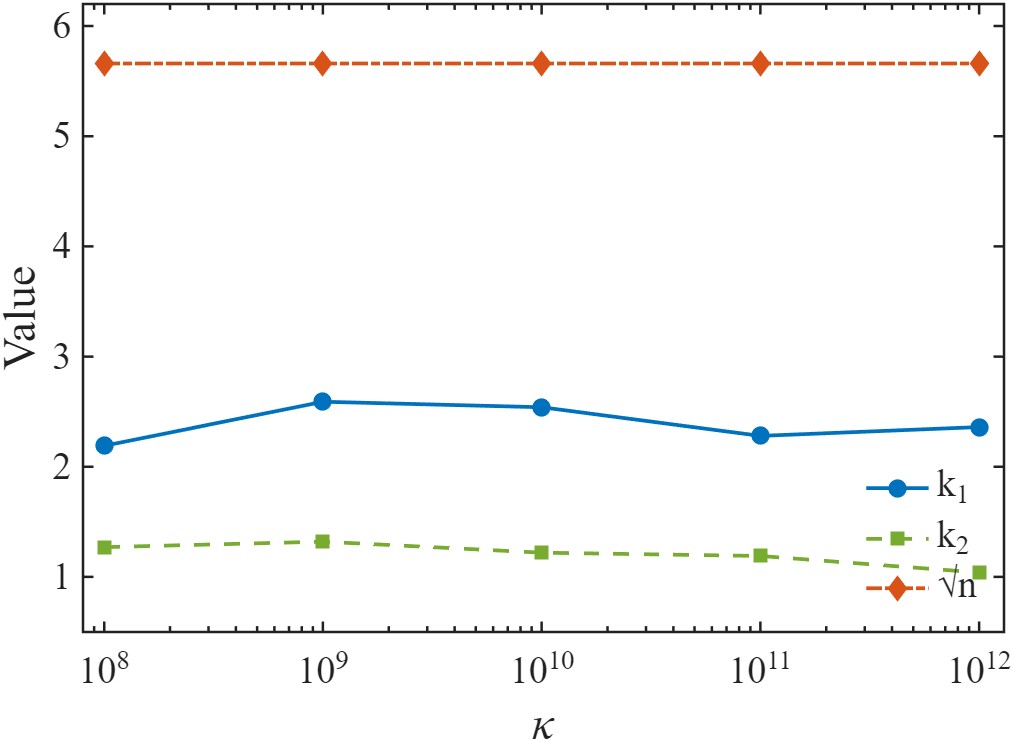}
\caption{$k_{1}, k_{2}$ with different $\kappa_{2}(X)$}
\label{fig:valuesk}
\end{minipage}
\end{figure}

\begin{figure}
\centering
\begin{minipage}{0.44\textwidth}
\centering
\includegraphics[width=\textwidth]{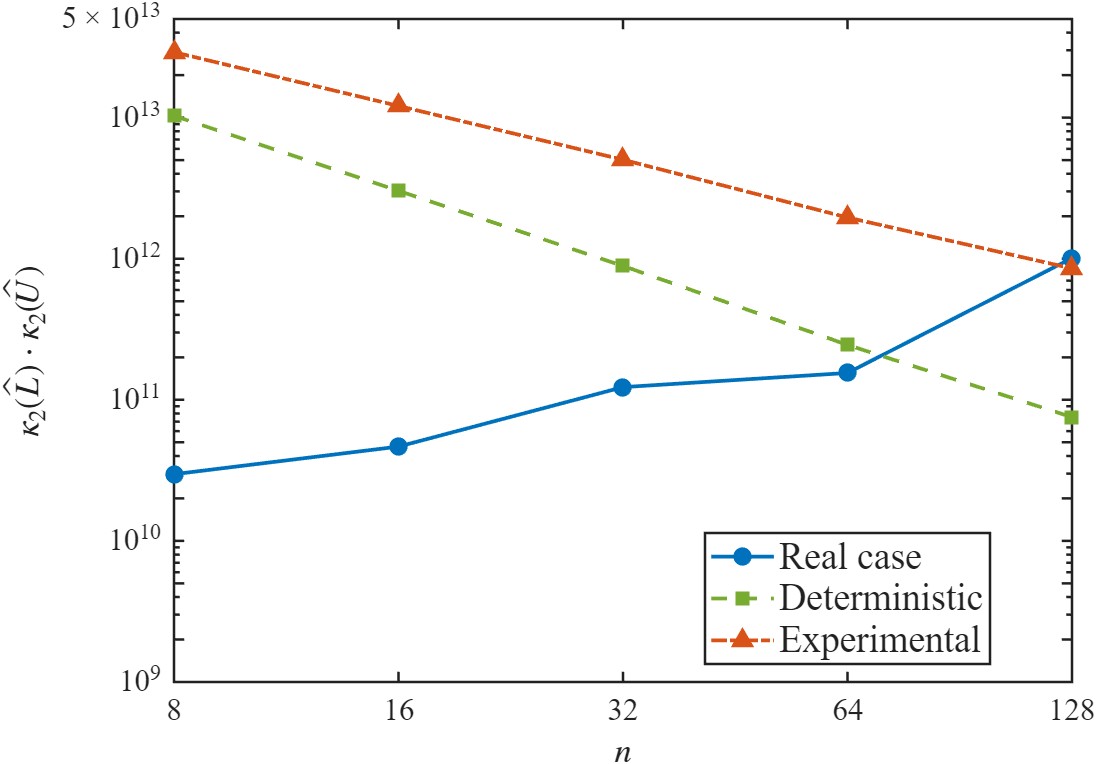}
\caption{Comparison of $\kappa_{2}(\hat{L})\kappa_{2}(\hat{U})$ with different $n$}
\label{fig:lun}
\end{minipage}
\hfill
\begin{minipage}{0.44\textwidth}
\centering
\includegraphics[width=\textwidth]{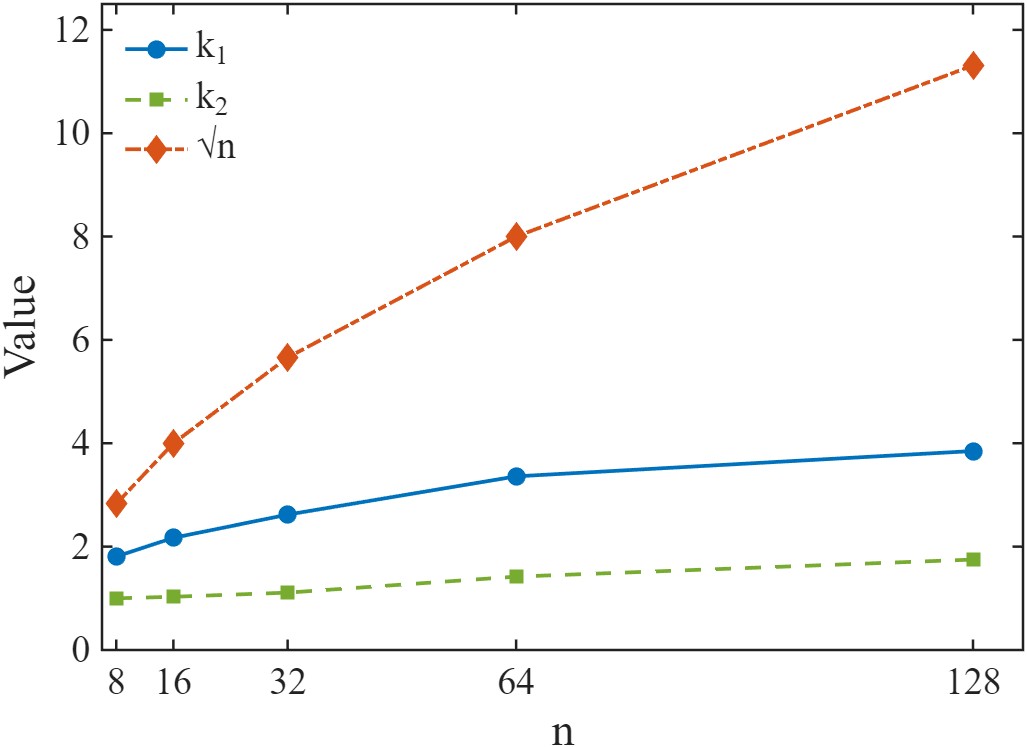}
\caption{$k_{1}, k_{2}$ with different $n$}
\label{fig:valuesn}
\end{minipage}
\end{figure}

In this section, we present numerical experiments to examine the sufficient condition for $\kappa_{2}(X)$ in \eqref{eq:s4}, as discussed in Section~\ref{sec:comparison}. We use $X \in \mathbb{R}^{m\times n}$ with $\epsilon=0.5$ and $p=0.6$, comparing the actual value of $\kappa_{2}(\hat{M})\kappa_{2}(\hat{N})$ with the deterministic bound $\frac{1}{21.4p_{2}}$ and the experimental bound $\frac{\sqrt{n}}{21.4p_{2}}$. We also compare $k_{1}$, $k_{2}$ and $\sqrt{n}$. To investigate the effect of $\kappa_{2}(X)$, we set $m=1032$, $n=32$ and $s=64$, while varying $\kappa_{2}(X)$. To study the effect of $n$, we fix $\kappa_{2}(X)=10^{10}$, $m=1024$ and $s=256$ while varying $n$. Numerical results are shown in Figure~\ref{fig:luk}-Figure~\ref{fig:valuesn}. $\kappa_{2}(\hat{L})\kappa_{2}(\hat{U})$, $\frac{1}{21.4p_{2}}$, and $\frac{\sqrt{n}}{21.4p_{2}}$ are referred to as 'Real case', 'Deterministic' and 'Experimental', respectively. Figures~\ref{fig:luk} and~\ref{fig:lun} verify that \eqref{eq:s4} is a sufficient condition. As $\kappa_{2}(X)$ or $n$ increases, the actual value of $\kappa_{2}(\hat{M})\kappa_{2}(\hat{N})$ eventually exceeds both the deterministic and experimental bounds. Moreover, the experimental bound applies to more practical cases than the deterministic bound. Figure~\ref{fig:valuesk} and Figure~\ref{fig:valuesn} show that $k_{1}$ and $k_{2}$ are substantially smaller than $\sqrt{n}$. This demonstrates that the sufficient condition in \eqref{eq:s4} is less restrictive than that in \cite{LUChol}.

\subsection{Tests of the bounds of residual for RCLUPP}
\begin{figure}
\centering
\begin{minipage}{0.48\textwidth}
\centering
\includegraphics[width=\textwidth]{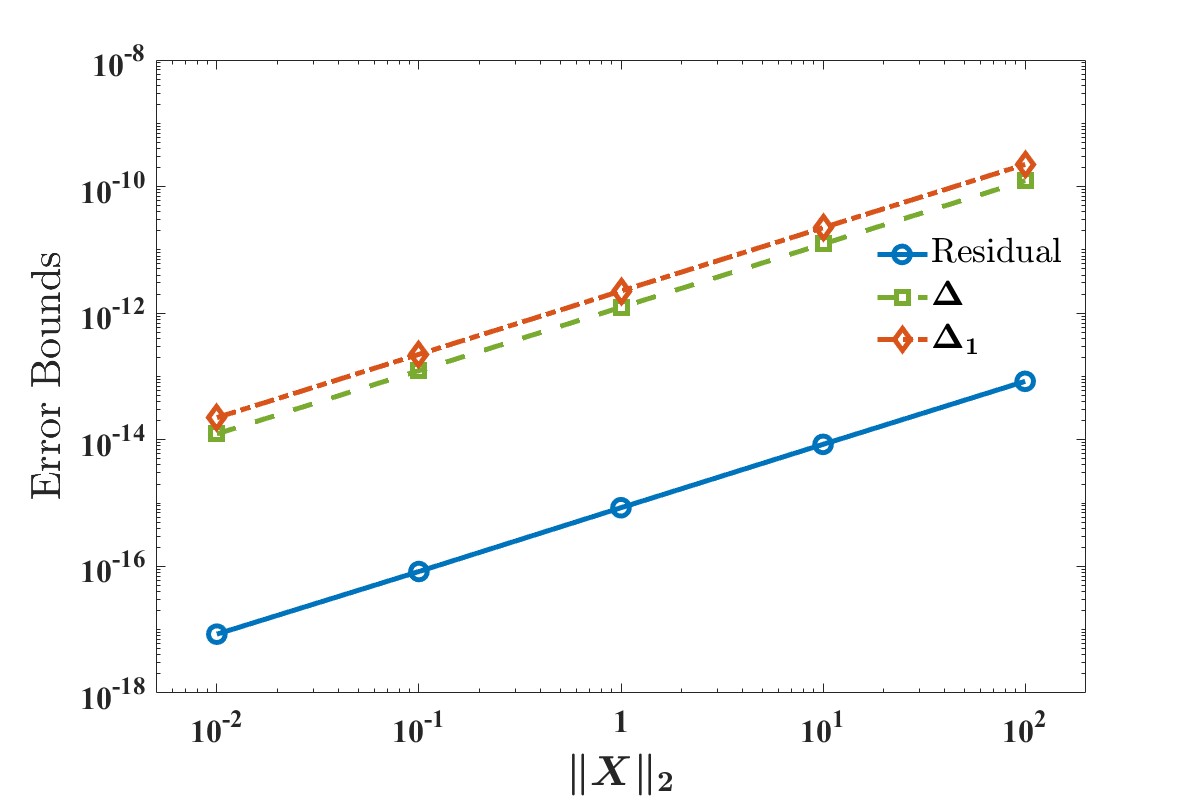}
\caption{Comparison of the bounds of residual with different $\kappa_{2}(X)$}
\label{fig:bx}
\end{minipage}
\hfill
\begin{minipage}{0.48\textwidth}
\centering
\includegraphics[width=\textwidth]{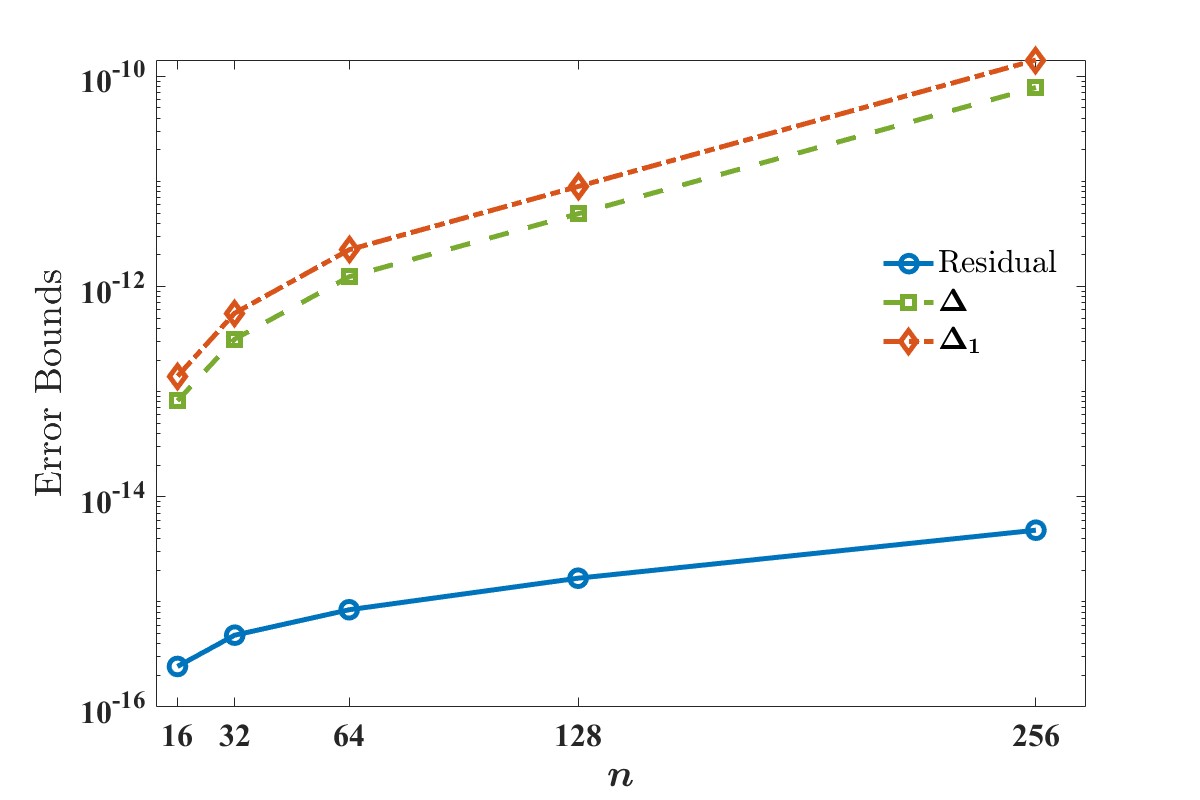}
\caption{Comparison of the bounds of residual with different $n$}
\label{fig:bn}
\end{minipage}
\end{figure}

To validate the sharper bound of residual $\Delta$ in \eqref{eq:res}, we perform a series of numerical experiments. When bounding the residual of RCLUPP following the approaches in the existing works, we replace $\norm{X}_{F}$ in $\Delta$ by $\sqrt{n}\norm{X}_{2}$, yielding $\Delta_{1}=\frac{7.38}{\sqrt{1-\epsilon}} \cdot hn^{2}\uu \cdot \norm{X}_{2}$ with $h=\frac{0.78}{\sqrt{1+\epsilon}}-\frac{0.1}{\sqrt{1-\epsilon}}$. We compare $\Delta$ and $\Delta_{1}$ with the actual computed residual to demonstrate the superiority of $\Delta$ over $\Delta_{1}$. Using the SVD-based matrix $X \in \mathbb{R}^{m\times n}$ with $\epsilon=0.5$, $p=0.6$, $\kappa_{2}(X)=10^{4}$ and $m=2048$, we investigate the effects of $\norm{X}_{2}$ and $n$. When varying $\norm{X}_{2}$, we fix $n=64$ and $s=128$. For different $n$, we set $\norm{X}_{2}=1$, $s=512$. The results are reported in Figure~\ref{fig:bx} and Figure~\ref{fig:bn}. It is clear that the sharper bound of residual $\Delta$ consistently provides a tighter estimate of the actual residual than $\Delta_{1}$, highlighting the advancement of the rounding error analysis in this work.

\subsection{Tests of the efficiency of the algorithms}
\begin{figure}
\centering
\begin{minipage}{0.48\textwidth}
\centering
\includegraphics[width=\textwidth]{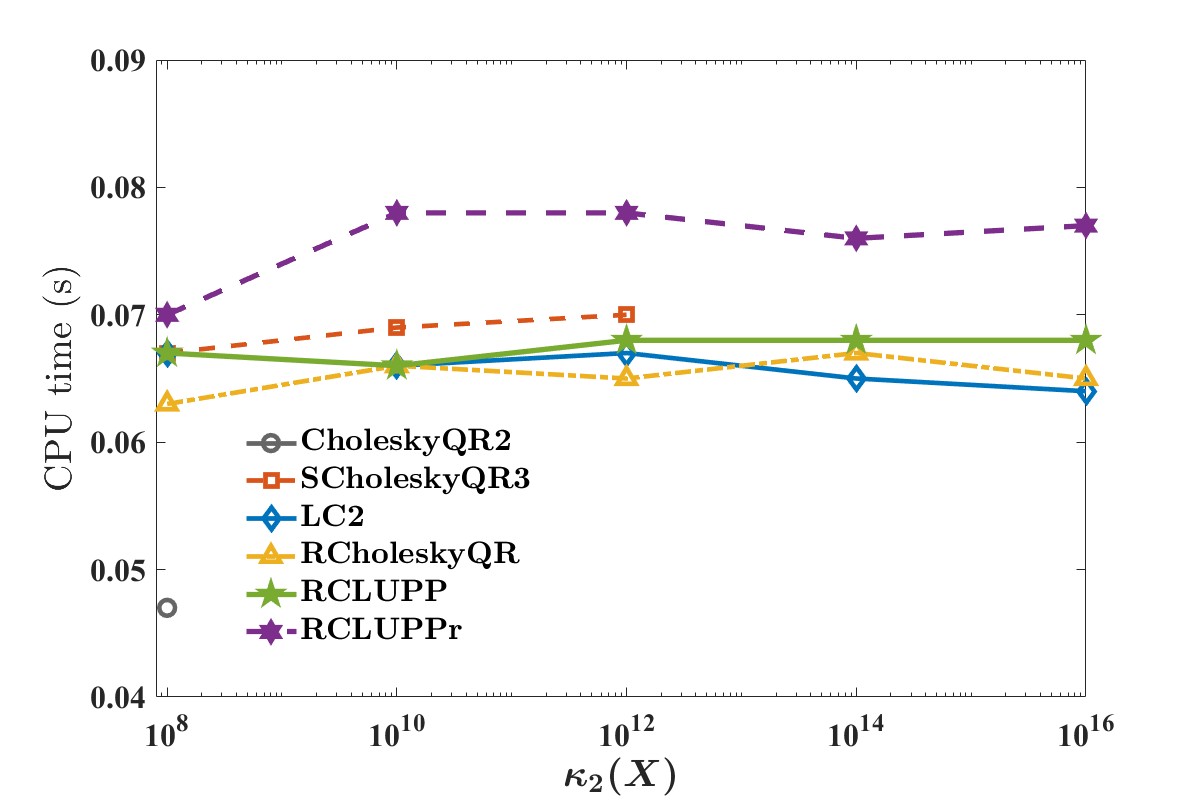}
\caption{Comparison of CPU time (s) with different $\kappa_{2}(X)$}
\label{fig:CPU0}
\end{minipage}
\hfill
\begin{minipage}{0.48\textwidth}
\centering
\includegraphics[width=\textwidth]{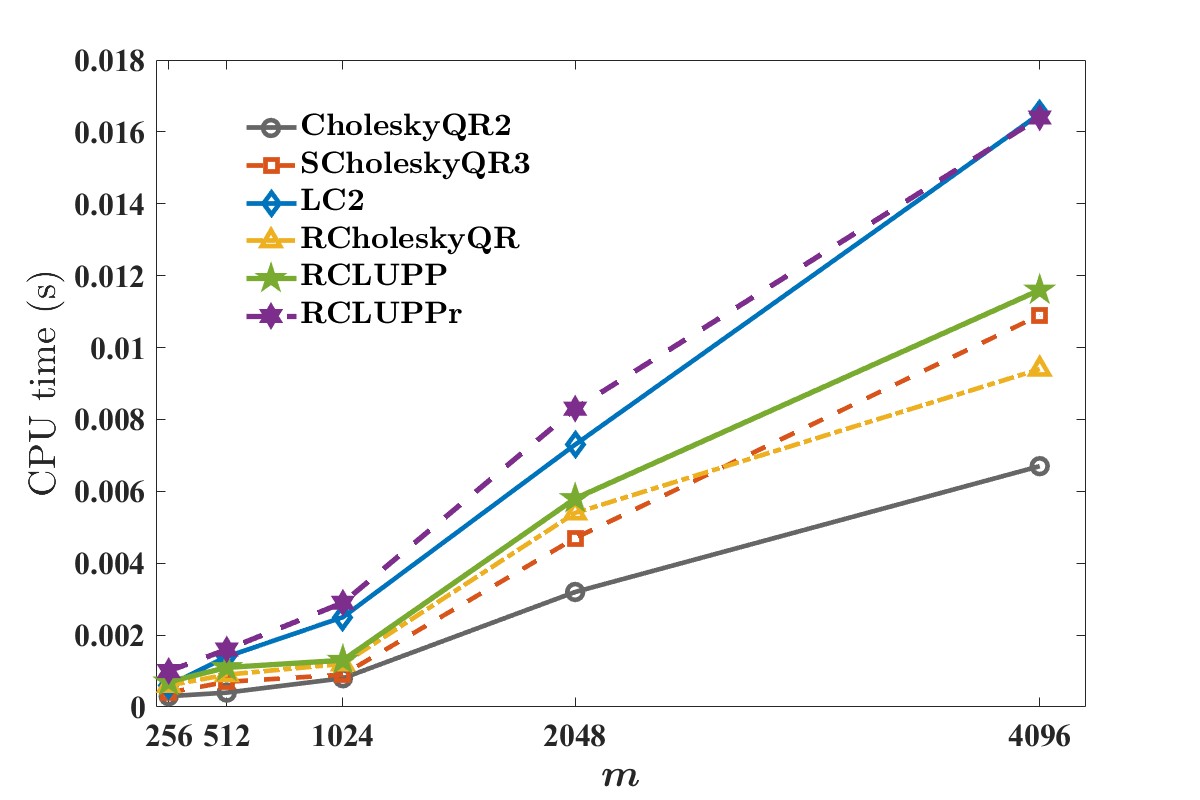}
\caption{Comparison of CPU time (s) with different $m$}
\label{fig:CPU1}
\end{minipage}
\end{figure}

\begin{figure}
\centering
\begin{minipage}{0.48\textwidth}
\centering
\includegraphics[width=\textwidth]{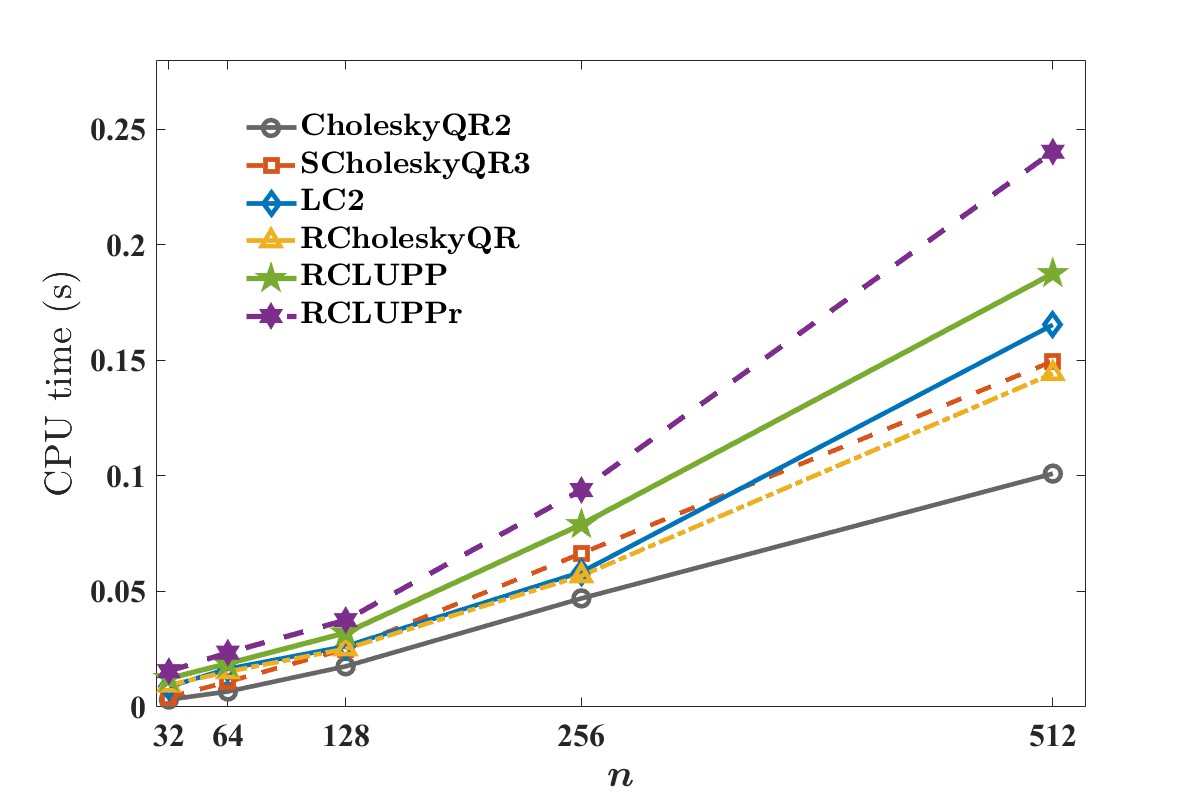}
\caption{Comparison of CPU time (s) with different $n$}
\label{fig:CPU2}
\end{minipage}
\hfill
\begin{minipage}{0.48\textwidth}
\centering
\includegraphics[width=\textwidth]{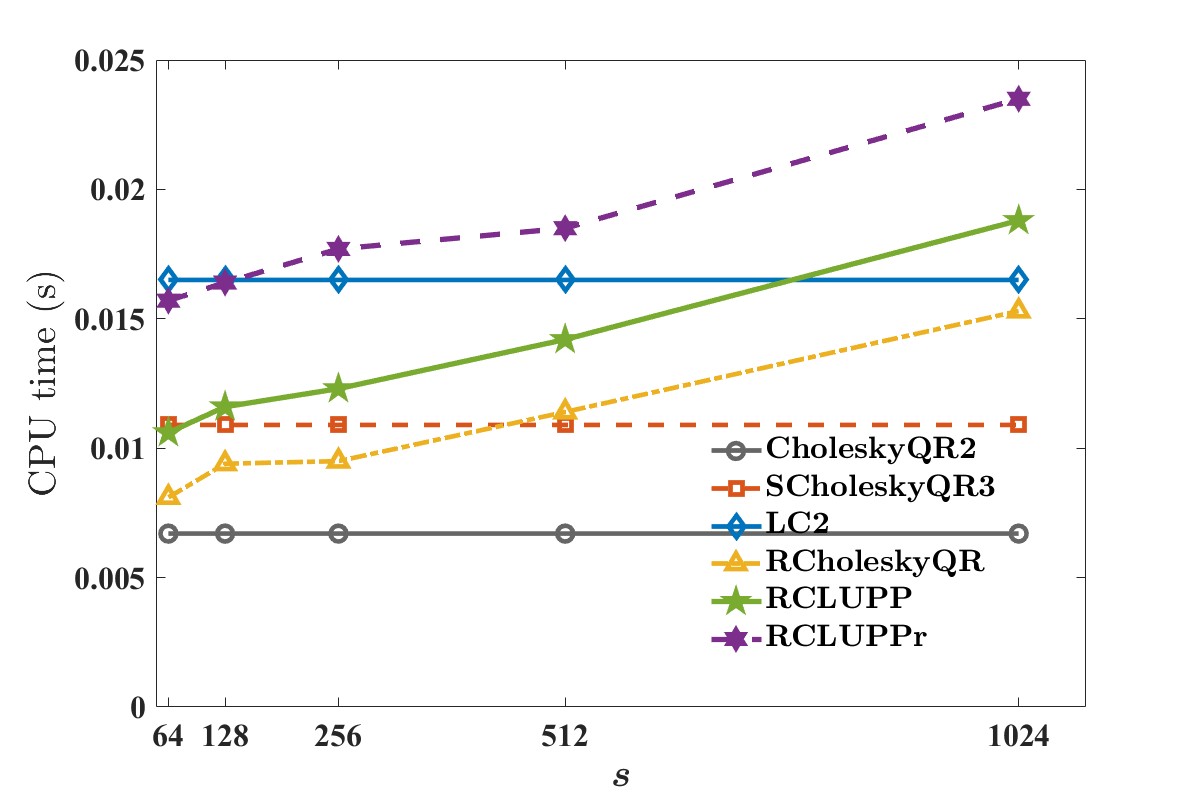}
\caption{Comparison of CPU time (s) with different $s$}
\label{fig:CPU3}
\end{minipage}
\end{figure}

In this section, we examine the computational efficiency of the algorithms through the CPU time (s). For RCLUPP, RCLUPPr and RCholeskyQR, we still take $\epsilon=0.5$ and $p=0.6$. Using the SVD-based matrices, we compare the CPU time (s) of RCLUPP and RCLUPPr against CholeskyQR2, SCholeskyQR3, LC2 and RCholeskyQR under varying $\kappa_{2}(X)$, $m$, $n$ and $s$ for the input $X \in \mathbb{R}^{m\times n}$. To study the effect of $\kappa_{2}(X)$, we fix $m=4096$, $n=256$ and $s=512$ for the randomized algorithms. When varying $m$, we set $\kappa_{2}(X)=10^{4}$, $n=64$ and $s=128$. For different $n$, we set $\kappa_{2}(X)=10^{4}$, $m=4096$ and $s=1024$. When varying $s$, we set $\kappa_{2}(X)=10^{4}$, $m=4096$ and $n=64$. Numerical results of the CPU time (s) are reported in Figure~\ref{fig:CPU0}-Figure~\ref{fig:CPU3}. Note that the time required to generate the sketch matrix $S$ is excluded from the reported CPU time (s) for the randomized algorithms.

Figure~\ref{fig:CPU0}-Figure~\ref{fig:CPU3} show that when $mn$ is large and $s$ is close to $n$, RCLUPP is more efficient than LC2. It is also faster than SCholeskyQR3 in some cases. Although RCLUPPr incurs the cost of LUP decomposition before matrix sketching, its CPU time (s) remains comparable to that of LC2 for the tall-skinny $X \in \mathbb{R}^{m\times n}$ with small $s$. Although RCLUPP and RCLUPPr are slightly slower than RCholeskyQR because of LUP decomposition, their efficiency is generally acceptable. In the real practice, setting $s=2n$ or even $n$ can further accelerate both algorithms without sacrificing the accuracy, which is a strategy in the implementation.

\section{Conclusions and discussions}
\label{sec:conclusions}
In this work, we introduce RCLUPP and its variant RCLUPPr, a new family of randomized CholeskyQR-type algorithms based on LUP decomposition. By strategically integrating matrix sketching, LUP decomposition and the thin HouseholderQR, our proposed RCLUPP and RCLUPPr achieve a significantly improved balance among the applicability to the ill-conditioned problems, accuracy, numerical stability and efficiency compared with the existing CholeskyQR-type algorithms. A detailed rounding error analysis of RCLUPP yields a sharper bound of residual than those in the literature. Extensive experiments on both synthetic and the real-world matrices validate the theoretical findings. The algorithms are shown to be particularly effective for QR decomposition of the tall-skinny matrices in the real applications.

The direction of our future work includes extensions of our algorithms to the rank-deficient cases, incorporation of the mixed-precision arithmetic for even greater performance on the modern hardware and the implementation of the algorithms on GPU. We believe that RCLUPP and RCLUPPr provide a practical and robust tool for QR factorization in the modern scientific computing.

\section*{Acknowledgement}
The authors declare that they have no conflict of interest. The work of Zhonghua Qiao is partially sponsored by NSFC/RGC Joint Research Scheme grant N\_PolyU5145/24 and Hong Kong Research Grants Council GRF grant 15305624.

\bibliography{references}

\end{document}